\documentclass[11pt, epsf]{article}
\fontfamily{times}
\usepackage{graphicx}
\usepackage{geometry, amsmath}
\usepackage{bm}
\include{epsfig}
\include{epsf}
\newtheorem{thm}{Theorem}
\newtheorem{cor}[thm]{Corollary}

\newtheorem{lemma}{Lemma}

\newcommand{\ben}{\begin{enumerate}}
\newcommand{\een}{\end{enumerate}}
\newcommand{\beq}{\begin{eqnarray}}
\newcommand{\eeq}{\end{eqnarray}}
\newcommand{\beqn}{\begin{eqnarray*}}
\newcommand{\eeqn}{\end{eqnarray*}}
\newcommand{\e}{\varepsilon}

\newcommand{\be}{\begin{equation}}
\newcommand{\ee}{\end{equation}}
\font\Bbb=msbm10

\def\sk1{\vskip 10pt}

\def\Bbb#1{I\!\! #1}

\def\I{\bold I}

\def\y{{\bold y}}

\def\x{{\bold x}}

\def\f{{\bold f}}
\def\e{{\bold e}}
\def\t{{\tt t}}
\def\u{{\bold u}}

\def\D{{\bold D}}

\def\b{{\bold b}}

\def\w{{\bold w}}

\def\S{{\bold S}}

\def\U{{\bold u}}

\def\b0{{\bold 0}}

\def\bold{\bf}

\def\t{{\bold t}}
\def\f{{\bold f}}
\def\bb{{\bold b}}

\def\BSigma{\boldsymbol{\Sigma}}
\def\thn{\boldsymbol{\hat\theta}_n}
\def\ths{\boldsymbol{\hat\theta}_n^*}
\def\thz{\boldsymbol{\theta}_0}

\def\tht{\boldsymbol{\theta}}

\def\ssum{\sum_{i=1}^{n}}

\def\thsts{\boldsymbol{\hat\theta}_{_{STS}}}
\def\thh{\boldsymbol{\hat\theta}}
\def\BSigma{\boldsymbol{\Sigma}}
\def\thn{\boldsymbol{{\hat\theta}_n}}
\def\ths{\boldsymbol{{\hat\theta}_n^*}}
\def\thz{\boldsymbol{\theta}_0}

\def\tht{\boldsymbol{\theta}}

\def\thnn{\hat\theta_{ni}}
\def\thss{\hat\theta_{ni}^*}
\def\thzz{\theta_0}
\def\thtt{\theta}
\def\jsum{\sum_{j=1}^{n_i}}
\def\isum{\sum_{i=1}^{N}}
\def\e{\epsilon}
\def\l{\left[}
\def\r{\right]}
\def\ll{\left(}
\def\rr{\right)}
\def\ssum{\sum_{i=1}^{n}}

\font\absmall=cmr12 at 11 pt
\def\Cal{\cal}

\font\absmall=cmr8 at 8 pt
\def\Cal{\cal}
\def\l{\left[}
\def\r{\right]}
\def\rr{\right)}
\def\ll{\left(}
\providecommand{\keywords}[1]{\textbf{\textit{Keywords: }} #1}

\newcommand{\qed}{\hspace*{\fill}Q.E.D.}  

\title{Recycled Two-Stage Estimation in Nonlinear Mixed Effects Regression Models}
\author{Ben Boukai\thanks{%
Email: bboukai@iupui.edu } \ and Yue Zhang\thanks{%
Email: yz65@umail.iu.edu}\\ Department of Mathematical Sciences, IUPUI\\
Indianapolis, Indiana, 46202 }

\begin{document}
\maketitle


\begin{abstract}
\noindent We   consider a re-sampling  scheme for estimation of the population parameters in the mixed effects nonlinear regression models of the type use for example in clinical pharmacokinetics, say. We provide an estimation procedure which {\it recycles}, via random weighting, the relevant two-stage parameters estimates to construct consistent estimates of the sampling distribution of the various estimates. We establish the asymptotic normality of the resampled estimates and demonstrate the applicability of the {\it recycling} approach in a small   
simulation study and via example. 
\end{abstract}


\bigskip

\setlength{\parindent}{0pt}

\section{Introduction} 
Hierarchical mixed-effects nonlinear regression models are widely used nowadays to analyze repeated measures observations. Data consisting of repeated measurements taken on each of a number of individuals arise commonly in biological and biomedical applications. Such models provide a natural settings for the analysis of data from pharmacokinetic studies obtained from a group of individuals which assumed to constitute a random sample from a relevant population of interest.

The hierarchical nonlinear model can be considered as an extension of the ordinary nonlinear regression models constructed to handle data obtained from several individuals. Modeling this kind of data usually involves a ``functional" relationship between at least one of the predictor variables, $x$, and the measured response, $y$,  within the individual's data.  As it often the case, the assumed 'functional' model between the response $y$ and the predictor $x$,  is based on some on physical or mechanistic grounds and is usually nonlinear in its parameters. These parameters are typically  estimated from the data by some techniques suitable for nonlinear regression. 

Figure \ref{1.1} below shows drug concentration by time  profiles for a study of the anti-asthmatic drug, {\it Theophylline}, as reported in Boeckmann, Sheiner and Beal (1994). Same dose of the drug was orally administered to 12 subjects, and over the subsequent 24 hour, serum concentrations were measured at ten time points per subject. For each subject, the pattern is of a rapid increase (post-drug) up to a to a peak concentration,  followed by an apparent exponential decay. A common pharmacokinetics model to describe this relation following an oral administration of the {\it Theophylline}  is the one-compartment model with first-order absorption and elimination rates (see for Example Davidian and Giltinan (1995)) .

\begin{figure}[h]
	\centering 
	\includegraphics[width=0.45 \textwidth]{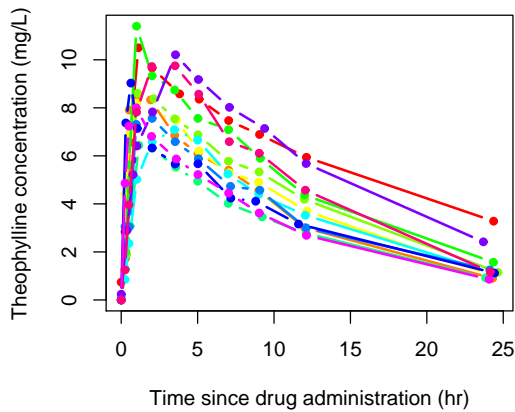}
	\caption{Drug concentrations for 12 participants in the {\it Theophylline} pharmacokinetics study }
	\label{1.1}
\end{figure}

As we can see from this figure, this type of data involves within-subject variability as well as between-subject variability from an assumed population pharmacokinetic model.  In such an hierarchical population model, fixed-effect parameters  quantify the population average kinetics of the drug whereas inter-individual random effect parameters quantify the magnitude of inter-individual variability. 

 The basic hierarchical linear regression model was pioneered by Sheiner, Rosenberg and Melmon (1972), which accounted for both types of variations; of within and between subjects. The nonlinear case received widespread attention in later developments. Lindstrom and Bates (1990) proposed a general nonlinear mixed effects model for repeated measures data and  proposed estimators combined least squares estimators and maximum likelihood estimators (under specific normality assumption). Vonesh and Carter (1992) discussed nonlinear mixed effects model for unbalanced repeated measures. Additional related references include: Mallet (1986), Davidian and Gallant (1993), Davidian and Giltinan (1993, 1995).

In all, the standard approach for statistical inference in hierarchical nonlinear models,  is typically based on full distributional assumptions for both,  the intra and inter individual random components. The most common assumption is that both random components are considered to be normally distributed. However, this can be a questionable assumption in many cases. Our main results in this work are built on more generalized assumptions in which the normally distributed random terms are not required.  

One of the main approaches for estimation in such hierarchical 'population' models  is the two-stage estimation methods. At the first stage to estimate the 'individual'-level parameters and then to combine them by some manner to obtain the 'population'-level parameters. However,  the main challenge to such two-stage estimation methods  is to obtain the sampling distributions of the final estimators in order to evaluate performance, especially when there is no sufficient data available or whenever existing asymptotic results are not accurate. For most part, the performance of these estimation methods can only be evaluated empirically, primarily via the so-called Monte-Carlo simulations-- see related references including: Sheiner and Beal (1981, 1982, 1983) and Davidian and Giltinan (1995, 2003). Hence, an alternate and more data oriented methodology should be considered. Bar-Lev and Boukai (2015) proposed a variant of the random weighting technique, which is termed herein {\it recycling}, as a valuable and valid alternative methodology for evaluation and comparison of the estimation procedure. Boukai and Zhang (2018) studied the asymptotic properties (asymptotic consistency and normality) of the {\it recycled} estimated in a one-layered nonlinear regression model.  

In this paper we extend to the hierarchical nonlinear regression models the approach of Bar-Lev and Boukai (2015) to include general random weights and with minimal (only moments) assumptions on the random error-terms/effects. In Section 2, we introduce and study the standard two-stage (STS) estimates in the hierarchical settings of nonlinear mixed effect models, and establish the asymptotic consistency and asymptotic normality of the STS estimators in such general settings. As far as we know, these are the first provably valid asymptotic distributional results concerning the STS estimation procedure in the context of hierarchical nonlinear regression. Furthermore, in Section 3 we introduce a specialized re-sampling scheme to obtain the {\it recycled} version of the STS estimators and demonstrate their  the asymptotic consistency and  normality as well. The results of extensive simulation studies and a couple of detailed illustrations are provided in Section 4. The proofs of our main results along with many other technical details are provided in Section 5.

\section{The Basic Hierarchical (Population) Model}
Consider a study involving a random sample of $N$ individuals, where the nonlinear regression model (as in Boukai and Zhang (2018)) is assumed to hold for each of the $i$-th individuals. That is, for each $i$, 
$i=1, 2, \dots, N$, we have available the $n_i$ (repeated) observations (with $n_i>p$) on the response variable in the form of $\y_i:=(y_{i1}, y_{i2}, \dots, y_{i n_i})^\t$, where
\be\label{4.1.1}
 y _{i j}=f(\x_{ij}; \, \tht_i)+\epsilon_{ij},   \ \ \ \     j=1,\dots, n_i, 
\ee
and  $ \x_{ij} $ is the $ j$-th fixed input (or condition) for the $i$-th individual, which gives rise to the response, $ y_{ij}$, for $j=1, \dots, n_i$ and $i=1, \dots, N$. Here, $f(\cdot)$ is a given nonlinear function and $\epsilon_{ij}$ denote some $i.i.d.$ $(0, \sigma^2)$ error-terms. That is, if we set ${\boldsymbol{\epsilon}}_{n_i}:=(\e_{i1},\e_{i2},\dots,\e_{in_i})^\t$, then 
$$
E({\boldsymbol{\epsilon}}_{n_i}) = \b0\ \ \ \text{and} \ \ Var({\boldsymbol{\epsilon}}_{n_i})\equiv Cov({\boldsymbol{\epsilon}}_{n_i} {\boldsymbol{\epsilon}}^\t_{n_i}) = \sigma^2\I_{n_i}.
$$
In the current context, the parameter vector $\tht=(\theta_1, \theta_2, \dots, \theta_p)^\t\in \Theta\subset {\Bbb R}$ can vary from individual to individual, so that $\tht_i$ is seen as the individual-specific realization of  $\tht$. More specifically, it is assumed that, independent of the error terms,  ${\boldsymbol{\epsilon}}_{n_i}$, 
$$
\tht_i:= \thz+\bb_i,
$$
where $\thz:=(\theta_{01}, \theta_{02}, \dots, \theta_{0p})^\t$,  is a fixed population parameter, though unknown, and $\bb_i=(b_{i1}, b_{i2}, \dots, b_{ip})^\t$ is a $p\times 1$ vector representing the random effects associated with $i$-th individual. It is assumed that the random effects, $\bb_1, \bb_2, \dots, \bb_N$ are independent and identically distributed random vectors satisfying,
$$
E(\bb_i)=\b0 \ \ \ \text{and} \ \ \ Var(\bb_i)\equiv Cov( \bb_i, \bb_i^t) = \D.
$$
Thus, $\tht_1, \tht_2, \dots, \tht_N$ are $i.i.d.$ random vectors with
$$
E(\tht_i)=\b0\ \ \  \text{and} \ \ \ Var(\tht_i) = \D. 
$$
 In the simple  hierarchical modeling it is often assumed that $\D$ is some diagonal matrix of the form $\D=Diag(\lambda_1^2,\lambda_2^2, \dots, \lambda_p^2)$ or even simpler, as $\D =\lambda^2 \I_p$ for some $\lambda>0$, and that both,   the error terms ${\boldsymbol{\epsilon}}_{n_i}$, and the random effects $\bb_i$ are normally distributed, so that, 
$$
{\boldsymbol{\epsilon}}_{n_i} \sim {\Cal N}_{n_i}(\b0, \sigma^2\I_{n_i}), \ \ \text{and} \ \ \ \bb_i\sim {\Cal N}_p(\b0, \D), 
$$
for each $=i=1, \dots, N$. In the more complex hierarchical modeling, more general structures of the {\it  within individual} variability $Var({\boldsymbol{\epsilon}}_{n_i})=\Gamma_i$ (for some $\Gamma_i$) and of the {\it between individuals} variability, $\D$, are possible.  However, even in the simplest structure, the available estimation methods for these model's parameters, $\thz, \sigma^2$ and $\D$ are typically  highly iterative in their nature and are based on the variations of the least squares estimation, and when available under some specific distributional assumptions, also on the maximum likelihood estimation procedures. In fact, many of the available results in the literature  hinge on the specific normality assumption and on the ability to effectively 'linearize' the regression function $f(\cdot)$ (see for example Bates and Watts (2007)). We point out that here {\bf we require no specific distributional assumptions (such as normality)} on either ${\boldsymbol{\epsilon}}_{n_i}$ nor $\bb_i$.
However, we focus attention on the {\it Standard Two Stage} (STS) estimation procedure advocated by Steimer, Golmard and Boisvieux (1984).

\section{The Two-Stage Estimation Procedure\label{2.0}}

For each $i=1, \dots, N$, let $\f_i(\tht)$ denote the $n_i\times 1$ vectors whose elements are $f(\x_{ij}, \tht), j=1, \dots, n_i$ then model (\ref{4.1.1}) can be written more succinctly as
\be
\y_i=\f_i(\tht_i)+{\boldsymbol{\epsilon}}_{n_i}
\ee
Accordingly, the STS estimation procedure can be described as follows: 
\begin{itemize}
  \item[]{\bf  On Stage I:}\ \  For each $i=1, \dots, N$ obtain $\hat\tht_{ni}$ as the minimizer of 
\be\label{4.3.3}
Q_i(\tht):= (\y_i- \f_i(\tht))(\y_i-\f_i(\tht))^\t\equiv \sum_{j=1}^{n_i}(y_{ij}-f(\x_{ij}, \tht))^2 ,
\ee
so as to form $\hat\tht_{n1}, \hat\tht_{n2}, \dots, \hat\tht_{nN}$, based on all the $M:=\sum_{i}^N n_i$ available observations. Next, estimate the {\it within-individual} variability component, $\sigma^2$, by
$$
\hat \sigma_M^2:=\frac{1}{M-pN} \sum_{i=1}^N Q_i({{\thn}_i}).
$$

  \item[]{\bf On Stage II:}\ \  Estimate the `population' parameter $\thz$ by
\be\label{4.2}
\thsts:=\frac{1}{N} \isum {\thn}_i. 
\ee
Next,  estimate $Var(\thsts)$ by $\S^2(\thh)/N$, where
$$
\S^2(\thh):= \isum({\thn}_i-\thsts)({\thn}_i-\thsts)^\t. 
$$
Finally estimate the {\it between-individual} variability component, $\D$, by 
\be\label{4.5}
\boldsymbol{ \hat{D} }=\S^2(\thh)- \min(\hat \nu ,  \hat \sigma_M^2) {\boldsymbol{\hat \Sigma}}_N, 
\ee
where ${\boldsymbol{\hat \Sigma}}_N:= \frac{1}{N}\isum {\boldsymbol{\Sigma}}_{n_i}({\thn}_i)$,  with ${\boldsymbol{\Sigma}}^{-1}_{n_i}$ defined as,
\be\label{b}
\BSigma^{-1}_n(\tht):= \frac{1}{n}\sum_{i=1}^{n}\nabla f_i(\tht)\nabla f_i(\tht)^\t, 
\ee
and where $\hat \nu$ is the smallest root of the equation $|\S^2_{_{STS}}-\nu {\boldsymbol{\hat \Sigma}}_N|=0$, see Davidian and Giltinan (2003) for details. 
 \end{itemize}

 Bar-Lev and Boukai (2015) provided a numerical study of this two-stage estimation procedure in the context of pharmacokinetics (hierarchical) modeling under the normality assumption. They also proposed a corresponding two-stage resampling (or recycling) algorithm, but based on ${\Cal{D}}irichlet(\boldsymbol{1})$ random weights. However, in this paper we consider a more general framework for the random weights to be used. 
 
 For each $n\geq 1$, we let the random weights, $\w_n=(w_{1:n}, w_{2:n}, \dots, w_{n:n})^\t$, be a vector of exchangeable nonnegative random variables with $E(w_{i:n})=1$ and   $Var(w_{i:n}):= \tau_n^2$,   and let $W_{i}\equiv W_{1:n}=(w_{i:n}-1)/\tau_n$ be the standardized version of $w_{i:n}$, $i=1, \dots, n$.  In addition we also assume, in similarity to Boukai and Zhang (2018) that,

\noindent {\bf \underline{Assumption W:}} The underlying distribution of the random weights $\w_n$ satisfies
\begin{enumerate}

  \item For all $n\geq 1$, the random weights $\w_n$ are independent of $(\epsilon_1, \epsilon_2, \dots, \epsilon_n)^\t$;
  \item $\tau^2_n=o(n)$, $E(W_iW_j)=O(n^{-1})$ and $E(W_i^2W_j^2)\to 1$   for all $i\ne j$,  $E(W_i^4)<\infty$ for all $i$.
\end{enumerate}
 
With such general random weights, the {\it recycled} version of the STS estimation procedure described in \ref{4.3.3}-\ref{b} above is:
\begin{itemize}
 \item[]{\bf  On Stage I$^*$:}\ \  For each $i=1, \dots, N$, independently generate random weights, $\w_i=(w_{i1},w_{i2},\dots,w_{in_i})^\t$ that satisfy {\it Assumption W} with $Var(w_{ij})=\tau^2_{n_i}$ and obtain ${\ths}_i$ as the minimizer of
\be\label{4.3.6}
Q_i^*(\tht):=  \sum_{j=1}^{n_i}w_{ij}(y_{ij}-f(\x_{ij}, \tht))^2 ,
\ee
so as to form ${{\ths}_1}, {\ths}_2, \dots, {\ths}_N$.

\item[]{\bf On Stage II$^*$:}\ \ Independent of {\bf Step I$^*$}, generate random weights,  $\u=(u_{1},u_{2},\dots,u_N)^\t$ that satisfy {\it Assumption W} with $Var(u_{i})=\tau^2_N$, and obtained the {\it recycled} version of $\thsts$  as:
\be\label{5.2}
\thsts^*:=\frac{1}{N} \isum u_i{\ths}_i
\ee
The {\it recycled} version $\D^*$ of $\D$ can be subsequently  obtained as described in {\bf Step II} above. 
\end{itemize}
\medskip\medskip

\section{Consistency of the Recycled STS Estimation Procedure}

In this section we present  some asymptotic results that establish and validate the consistency of the {\it recycled} STS estimator for general random weights satisfying the premises of {\it Assumption W}. We establish there results without the 'typical' normality assumption on the {\it within-individual} error terms, $\epsilon_{ij}$, nor on the {\it between-individual} random effects $\bb_i$. However, for simplicity of the exposition, we state these results in the case of $p=1$, so that $\Theta\in {\Bbb R}$. With that in mind, we denote for each $i=1, \dots, N$,   
$$
f_{ij}(\theta)\equiv f(x_{ij}, \theta), \ \ \text{for} \ \ j=1, \dots, n_i.
$$
Accordingly, the least squares criterion in (\ref{4.3.1}), becomes
 \[
Q_{ni}(\thtt):= \jsum(y_{ij}-f_{ij}(\theta))^2, 
\] 
and the LS estimator $\hat{\theta}_{ni}$ is readily seen as the solution of 
\be\label{4.1}
Q_{ni}'(\thtt):= 2\jsum\phi_{ij}(\thtt)=0
\ee
where,   
\be\label{4.1.1.1}
\phi_{ij}(\thtt):=-(y_{ij}-f_{ij}(\thtt))f^{'}_{ij}(\thtt), \ \   \ \ 
\ee
with $f^{'}_{ij}(\thtt):= d f_{ij}(\theta)/d \theta$, for $j=1\dots, n_i$ and for each $i=1\dots, N$. We write $f^{''}_{ij}(\thtt):= d f^\prime_{ij}(\theta)/d \theta$ and 
$\phi^\prime_{ij}(\thtt):= d \phi_{ij}(\theta)/d \theta$, etc.  As in Boukai and Zhang (2018), we also assume that $ f_{ij}^{'}(\thtt)$ and $f_{ij}^{''}(\thtt)$ exist for all $ \thtt $ near $\thzz $.  
However, to account for the inclusion  of the $(0, \lambda^2)$ random effect term, $b_i$,  in the model, we also  assume that, 

\noindent {\bf \underline{Assumption A: }} For each $i=1, \dots, N$
	\begin{enumerate}
		\item $a_{n_i}^{2}:=\sigma^2\jsum E(f^{'2}_{ij}(\thzz+b_i))\to \infty\ \ as \ \ n_i\to \infty$, ; 
		\item  $\underset{n_i\to\infty}{\limsup} \ \ a_{n_i}^{-2}\jsum \underset{|\thtt-\thzz-b_i|\le\delta}{\sup} f_{ij}^{''2}(\theta) <\infty $
		
		\item  $ a_{n_i}^{-2}\jsum f_{ij}^{'2}(\theta) \to \frac{1}{\sigma^2}$	uniformly in $|\theta-\thzz-b_i|\le \delta$.
	\end{enumerate}

\noindent In the following two Theorems we establish, under the conditions of {\it Assumption A}, the asymptotic consistency and  normality of $\hat\theta_{_{STS}}$.  Their proofs and some related technical results are given in Section \ref{5.33} below. 

\begin{thm}\label{4.3.1}
	Suppose that {\it Assumption A} holds, then there exists a sequence $\thnn$ of solutions of (\ref{4.1}) such that 
	\[
	\thnn=\thzz+b_i+a_{ni}^{-1}T_{ni}
	\]
	where $|T_{ni}|<K$ in probability, for each $i=1,2,\dots,N$. Further, there exists a sequence $\hat\theta_{_{STS}}$ as expressed in (\ref{4.2}) such that 
	\[
	\hat\theta_{_{STS}}-\thzz\overset{p}{\to}0, 
	\]
as $n_i\to\infty$, for $i=1,2,\dots,N$, and as $N\to \infty$. 
\end{thm}

\begin{thm}\label{4.3.2} Suppose that {\it Assumption A} holds. If
$$
\underset{N,ni\to \infty}{\lim} N/a_{ni}^2<\infty, 
$$
for all $i=1,2,\dots,N$, then there exists a sequence $\hat\theta_{_{STS}}$  as expressed in (\ref{4.2}) such that
	\[
	\hat\theta_{_{STS}}-\thzz= \frac{1}{N}\isum b_i -\psi_{_N,n_i },
	\]
	where $\sqrt{N}\psi_{_N,n_i}\overset{p}{\to}0$ . Further, 
	\[
	{\Cal{R}_N}:=\frac{\sqrt{N}}{\lambda}(\hat\theta_{_{STS}}-\thzz)\Rightarrow {\cal N}(0,1) 
	\]
	as $n_i\to\infty$, for  $i=1,2,\dots,N$, and as $N\to \infty$.  
\end{thm}

For the {\it recycled} STS estimation procedure as described in Section \ref{2.0} above, the {\it recycled}  version $\hat{\theta}^*_{ni}$ of $\hat{\theta}_{ni}$ is the minimizer of (\ref{4.3.6}), or alternatively, the direct solution of 
\be\label{4.8}
Q_{i}^{*\prime}(\thtt):= 2\jsum w_{ij}\phi_{ij}(\thtt)=0, 
\ee
where $\w_i=(w_{i1},w_{i2},\dots,w_{in_i})^\t$ are  the randomly drawn weights (satisfying {\it Assumption W}), for the $i$th individual, $i=1, 2, \dots, N$.   For establishing  comparable results to those given in Theorems \ref{4.3.1} and \ref{4.3.2} for the {\it recycled} version, 
 $\hat\theta^*_{_{STS}} =\isum u_i \hat{\theta}^*_{ni}/N$ of $\hat\theta_{_{STS}}=\isum \hat{\theta}_{ni}/N$, with the random weights $\u=(u_1, u_2, \dots, u_N)^\t$ as in {\bf Stage II$^*$}, we need the following additional assumptions.

\noindent {\bf \underline{Assumption B: }} In addition to {\it Assumption A},  we assume that $E(\e_{ij}^4)<\infty$ and that for each $i=1, 2, \dots,N$, 
	\begin{enumerate}
		\item $\underset{n_i\to\infty}{\limsup} \ \ a_{n_i}^{-2}\jsum \underset{|\thtt-\thzz-b_i|\le\delta}{\sup} f_{ij}^{'4}(\theta) <\infty$, 
		\item $\underset{n_i\to\infty}{\limsup} \ \ a_{n_i}^{-2}\jsum \underset{|\thtt-\thzz-b_i|\le\delta}{\sup} f_{ij}^{''4}(\theta) <\infty$, 
		\item As $n_i\to \infty$,  \ \ ${n_i}{a_{n_i}^{-2}}\to c_i\ge 0$ .
	\end{enumerate}

\noindent In Theorems \ref{5.1.1} and \ref{5.1.2} below we establish, under the conditions of {\it Assumptions A and B}, the asymptotic consistency and  normality of the {\it recycled} estimator $\hat\theta^*_{_{STS}}$.  Their proofs and some related technical results are given in Section \ref{5.44} below. 

\begin{thm}\label{5.1.1}
	Suppose that {\it Assumptions A} and {\it B} hold. Then there exists a sequence $\hat\theta_{ni}^*$ as the solution of (\ref{4.8}) such that
	\[
	\thss=\thnn+a_{ni}^{-1}T^*_{ni}
	\]
	where $|T_{ni}^*|<K\tau_{n_i}$ in probability, for $i= 1,\dots, N$. Further for any $\epsilon>0$, we have 
	\[
	P^*(|\hat\theta^*_{_{STS}}-\thzz|>\e)=o_p(1),
	\]
as $n_i\to\infty$, for $i=1,2,\dots,N$, and as $N\to \infty$.
\end{thm}

\begin{thm}\label{5.1.2}  Suppose that {\it Assumptions A} and {\it B} hold. If  for each $i=1, 2, \dots, N$, 
	$$
\frac{\tau_{n_i}}{\tau_{_N}}=o(\sqrt{n_i}), 
$$
then we have
	\[
	\hat\theta^*_{_{STS}}-\hat\theta_{_{STS}}= \frac{1}{N}\isum (u_i-1)\thnn -\psi^*_{_N,n_i}, 
	\]
	where $\frac{\sqrt{N}}{\tau_{_N}}\psi^*_{_N,n_i}\overset{p^*}{\to} 0$ as $N,n_i\to\infty$.  Additionally,  
	\[
	{\Cal{R}^*_N}:=\frac{\sqrt{N}}{\lambda \tau_{_N}} (\hat\theta^*_{_{STS}}-\hat\theta_{_{STS}})\Rightarrow {\cal N}(0,1),  
	\]
	as $n_i\to\infty$, for $i=1,2,\dots,N$, and as $N\to \infty$.  
\end{thm}
\noindent The proofs of Theorems \ref{5.1.1} and \ref{5.1.2} and some related technical results are given in Section \ref{5.44} below. The following corollary is an immediate consequence of the above results. It suggest that the sampling distribution of $\hat\theta_{_{STS}}$ can be well approximated by that of the {\it recycled} or re-sampled version of it, $\hat\theta^*_{_{STS}}$.

\begin{cor}
For all $t\in \Bbb{R}$, let   
$$
{\Cal{H}}_N(t)=P\left({\Cal{R}_N}\leq t\right),    \ \ \ \text{and} \ \ \ \ {\Cal{H}}_N^*(t)=P^*\left({\Cal{R}^*_N}\leq t\right),
$$
denote the corresponding c.d.f of ${\Cal{R}_N}$ and ${\Cal{R}^*_N}$, respectively. Then by Theorems \ref{4.3.2} and  \ref{5.1.2}, 
$$
	\underset{t}{\sup}| {\Cal{H}}_n^*(t)-{\Cal{H}}_n(t)|\to 0 \ \ \ in \ \ probability.
	$$
\end{cor}

\section{Implementation and Numerical Results} 
\subsection{Illustrating the STS Estimation Procedure}

To illustrate the main results of Section 4 for the hierarchical nonlinear regression model and the corresponding STS estimation procedure as described in \ref{4.3.3}-\ref{b} above, we consider a typical compartmental modeling from pharmacokinetics. In characterizing the pharmacokinetics of a drug disposition in the body, it is common to represent the body as a system of compartments and to assume that rates of transfer between compartments follow first-order or linear kinetics. Standard solution of the resulting differential equations shows that the relationship between drug concentration, as measured in the plasma and time (since administration of the drug to the body) may be described by a sum of exponential terms. For the standard two-compartment model, this relationship between the measure drug concentration $C(t)$ and the post-dosage time $t$, (following an intravenous administration), can be described through the nonlinear function of the form:
$$
f(t; {\boldsymbol{\eta}})=Ae^{-\alpha t}+Be^{-\beta t},
$$ 
with ${\boldsymbol{\eta}}:=(A, \alpha, B, \beta)^\prime$ is a parameter representing the various kinetics rate constants, such as the rate of elimination, rate of absorption, clearance, volume, etc. Since these constants (i.e. parameters) must be positive, we re-parametrize the model with $\tht\equiv \log(\boldsymbol{\eta})$,  so that with $t>0$,   
\be\label{6.1.1.1}
f(t; \tht)=exp(\theta_{1})exp\{-exp(\theta_{2})t\}+exp(\theta_{3})exp\{-exp(\theta_{4})t\}, 
\ee
with $\tht=(\theta_1, \theta_2, \theta_3, \theta_4)^\t\in {\Bbb R}^4$.  For the simulation stdy we conducted here, we consider a situation in which the (plasma) drug concentrations $\{y_{ij} \}$ of $N$ individuals were measure at post-dose times $t_{ij}$ and are related as in
model (\ref{4.1.1}) via the nonlinear regression model,
$$
y_{ij}=f(t_{ij}; \tht_{i})+\epsilon_{ij}, 
$$
for $j=1, \dots, n_i$ and $i=1, \dots, N$. Here, as in Section 4, $\epsilon_{ij}$ are the standard $(0, \sigma^2)$ error terms and $\tht_i=\thz+\bb_i$, where $\bold{b}_i$ are independent identically distributed random effects terms,  with mean $\bold{0}$ and unknown variance $\lambda^2{\boldsymbol{I}}_{4\times 4}$.  Accordingly, we have in all a total of 6 unknown parameters, namely, $\thz=(\theta_{10}, \theta_{20}, \theta_{30}, \theta_{40})^\t, \ \sigma$ and $\lambda$.  

 Since $\sigma$ and $\lambda$ represent variation within and between individuals (respectively), different setting for these two lead to very different situations.  For instance, Figure 1(a) below, depicts the situation for $N=5$ and $n_i\equiv n =15$, each,  when $\sigma=0.1$ and $\lambda=0.1$, so that the variation between individuals are similar to variation within individuals. Figure 1(b) depicts the situation with $\sigma=0.05, \lambda=1$, so that the variation between individuals is much larger than variation within individuals. 
 
 For the simulation, we set $\thz=(1,0.8,-0.5,-1)^\t$,  and for each $i$, the times $t_{ij}, j=1,\dots,n$ were generated uniformly from $[0,8]$ interval. To allow for different 'distributions', the error terms, $\e_{ij}$, as well as the random effect terms,  $\bb_i$, were generated either from the {\bf (a)} {\it Truncated Normal}, {\bf (b)} {\it Normal} and {\bf (c)} {\it Laplace} distributions -- all in consideration of {\it Assumption A} in our main results.
 
 For each simulation run, with the {\it Truncated Normal} distribution for the error-terms and the random effects terms,   we calculated the value of  $\hat\tht^k_{_{STS}}$ as an estimator of $\thz$ and repeated this procedure $M=1,000$ times to calculate the corresponding Mean Square Error (MSE) as followed,
\[
MSE=\frac{1}{M}\sum_{k=1}^{M}||\hat \tht^k_{_{STS}}-\thz||^2
\]
\vfill\eject

The corresponding simulation results obtained for various values of $N$ and $n$,  are presented in Table \ref{5.10} for $\sigma=0.1, \lambda=0.1$ and in  Table \ref{5.11} for $\sigma=0.05, \lambda=1$. 
 
\begin{figure}[h]
\begin{center}
\centerline{{\raise 1.0in \hbox{\epsfysize=2.6in \epsfbox{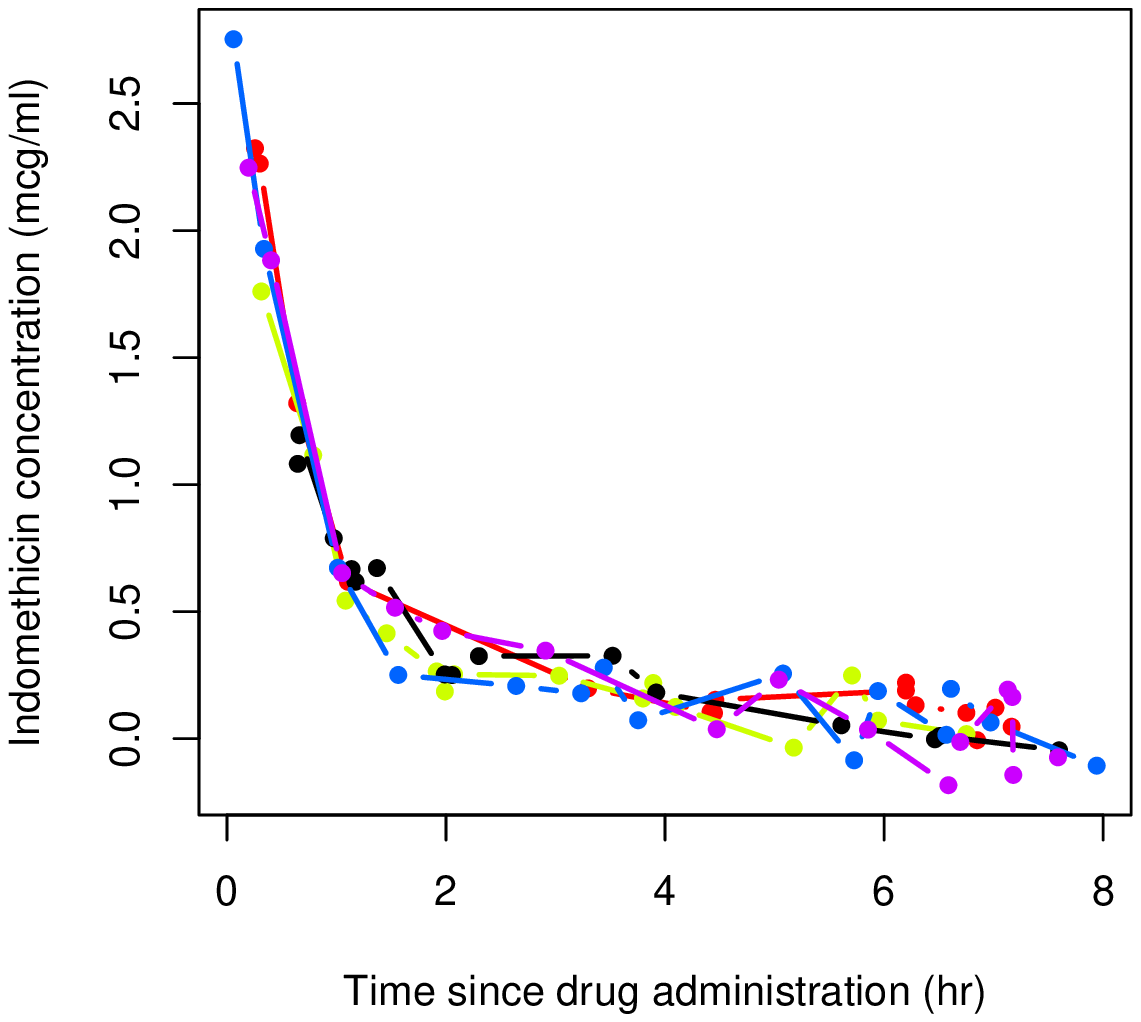}}}
{\raise 1in \hbox{\epsfysize=2.6in  \epsfbox{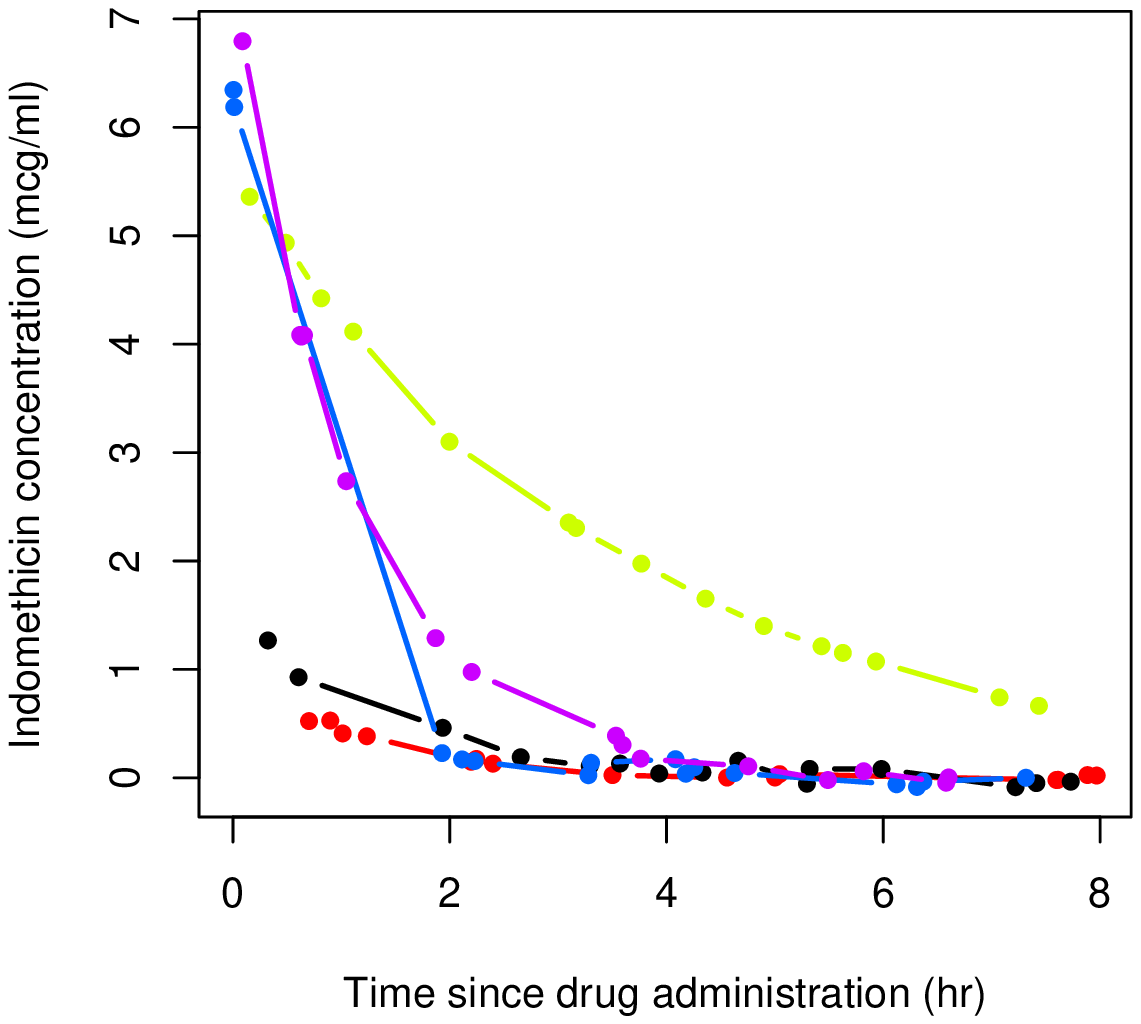}}}} \vskip-60pt {%
\caption{Drug plasma concentration vs time for (a) $\sigma=0.1, \lambda=0.1$; and for (b)  $\sigma=0.05, \lambda=1$ }}
\end{center}
\end{figure}


\begin{table}[!h]
\begin{center}
	{\absmall
		\begin{tabular}{c c c c c c}
			\hline\hline
			& n=15 &n=30 &n=50 &n=100 &n=200\\		
			\hline
			N=15     & 0.86616 &0.22885 &0.04651 &0.01141 &0.00632 \\
			\hline
			N=30     & 0.57666 &0.10713 &0.02442 &0.00573 &0.00334\\ 
			\hline
			N=50     & 0.45840 &0.08933 &0.02097 &0.00383 &0.00195\\
			\hline
			N=100    & 0.37852 &0.06918 &0.01245 &0.00216 &0.00103\\
            \hline
            N=200   & 0.35059 &0.05904 &0.00891 &0.00143 &0.00058\\
			\hline\hline		
				
		\end{tabular}
		}	
		\caption{The MSE of STS estimates for {\it truncated Normal} error-terms/effects with $\sigma=0.1, \lambda=0.1$.}
		\label{5.10}
\end{center}
\end{table}
From these two table, we see that with $n$ and $N$ both increasing, the MSE is decreasing, as expected. However, $\sigma=0.05, \lambda=1$ as in Table \ref{5.11}, $n$ increasing for a fixed $N$, doesn't contribute to smaller MSE, which is consistent with our main result Theorem \ref{4.3.1}, the STS estimate is not consistent with only $n_i\to \infty$, (this effect is more obvious in the case $\lambda$ is relatively large, as in the case of Table \ref{5.11}).

\begin{table}[!h]
\begin{center}
{\absmall
		\begin{tabular}{c c c c c c}
			\hline\hline
			& n=15 &n=30 &n=50 &n=100 &n=200\\		
			\hline
			N=15     & 1.00012 &0.63825 &0.56880 &0.47304 &0.46024\\
			\hline
			N=30     &0.69974 &0.39503 &0.33145 &0.35228 &0.32632\\ 
			\hline
			N=50     & 0.55675 &0.29437 &0.25938 &0.25004 &0.23474\\
			\hline
			N=100    &0.39821 &0.22447 &0.20213 &0.19734 &0.21995\\
			\hline
			N=200    &0.34921 &0.19447 &0.17476 &0.18824 &0.19581 \\
			\hline\hline			
		\end{tabular}
		}
		\caption{The MSE of STS estimates for {\it truncated Normal} error-terms/effects with $\sigma=0.05, \lambda=1$.}
\label{5.11}
\end{center}
\end{table}

For simulating the results of Theorem \ref{4.3.2}, we choose $\theta_{2}$ to be the unknown parameter, and use the main result to construct $95\%$ Confidence Interval as
\[
(\hat{\theta}_{_{STS}}-1.96\frac{\hat{\lambda}}{\sqrt{N}},\ \hat{\theta}_{_{STS}}+1.96\frac{\hat{\lambda}}{\sqrt{N}}),
\]
where 
\[
\hat{\lambda}^2=\frac{1}{N-1}\isum (\thnn-\hat\theta_{_{STS}})^2.
\]

The estimate for $\hat{\lambda}$ used here is the simple STS estimate, not the corrected one as in (\ref{4.5}). M=1,000 replications of such simulations were executed to determine the percentage of times the true value of the parameter estimates was contained in the interval.  We use $\sigma=0.5, \lambda=0.5$ and observed Coverage Percentages are provided in Table \ref{5.12} below.

\begin{table}[!h]
{\absmall	
\begin{center}
		\begin{tabular}{c c c c c c}
			\hline\hline
			& n=15 &n=30 &n=50 &n=100 &n=200\\		
			\hline
			N=15     &  0.903 &0.934 &0.933 &0.931 &0.931 \\
			\hline
			N=30     &0.896 &0.940 &0.940 &0.943 &0.944\\ 
			\hline
			N=50     &0.883 &0.941 &0.959 &0.944 &0.944\\
			\hline
			N=100 &0.828 &0.948 &0.946 &0.941 &0.944 \\
			\hline
			N=200 &0.759 &0.943 &0.932 &0.935 &0.949 \\
			\hline\hline			
		\end{tabular}
		\caption{Coverage Percentage of the CI for the {\it truncated Normal}  error-terms/effects with $\sigma=0.5, \lambda=0.5$.}
		\label{5.12}
\end{center}
}
\end{table}
From these results we can observe that with $n$ and $N$ both increase, the Coverage Percentage approximate to 0.95. While when $n$ is small (15), with $N$ increase, the Coverage Percentage is drifting farther away from the desired level of 0.95. This finding is consistent with our main result, the convergence require the condition $\underset{N,ni\to \infty}{\lim} N/a_{ni}^2<\infty$, which in this case becomes $\underset{n\to \infty}{\lim}\frac{1}{n}a_{n}^2/\sigma^2<\infty$, that is $\underset{N,n\to \infty}{\lim} N/n<\infty$ is required. Hence, when $N$ is much large than $n$, this condition does not hold.  Although for this model, error terms that follow the normal distribution do not satisfy {\it  Assumption A},  we used normal error terms in the simulations, and reported the resulting MSE and Coverage Percentage for 95\% confidence interval is in Table \ref{5.14} and Table \ref{5.15}. From the results we can observe that with $n$ and $N$ increasing,  the MSE are smaller and Coverage Percentage are closer to 0.95.
\begin{table}[!h]
{\absmall
\begin{center}
		\begin{tabular}{c c c c c c}
			\hline\hline
			& n=15 &n=30 &n=50 &n=100 &n=200\\		
			\hline
			N=15     &0.77176 &0.17458 &0.07880 &0.01116 &0.00615 \\
			\hline
			N=30     & 0.55483 &0.11852 &0.02966 &0.00605 &0.00324\\ 
			\hline
			N=50     & 0.47721 &0.09277 &0.02164 &0.00437 &0.00195\\
			\hline
			N=100    & 0.38275 &0.07416 &0.01217 &0.00231 &0.00104\\
			\hline
			N=200   & 0.33843 &0.05627 &0.00892 &0.00140 &0.00059\\
			\hline\hline			
		\end{tabular}
			\caption{The MSE of STS estimates for {\it Normal} error-terms/effects with $\sigma=0.1, \lambda=0.1$.}
			\label{5.14}
\end{center}
}
\end{table}

\begin{table}[!h]
{\absmall
\begin{center}
		\begin{tabular}{c c c c c c}
			\hline\hline
			& n=15 &n=30 &n=50 &n=100 &n=200\\		
			\hline
			N=15     &  0.918 &0.927 &0.939 &0.951 &0.922 \\
			\hline
			N=30     &0.901 &0.939 &0.944 &0.931 &0.932\\ 
			\hline
			N=50     &0.871 &0.947 &0.949 &0.950 &0.944\\
			\hline
			N=100  &0.851 &0.950 &0.934 &0.949 &0.948 \\
			\hline 
			N=200   & 0.740 &0.949 &0.944 &0.951 &0.945 \\
			\hline\hline			
		\end{tabular}
		\caption{Coverage Percentage of the CI for the {\it Normal}  error-terms/effects with $\sigma=0.5, \lambda=0.5$.}
		\label{5.15}
\end{center}
}
\end{table}

\begin{figure}[!h]
	\centering 
	\includegraphics[width=0.85\textwidth]{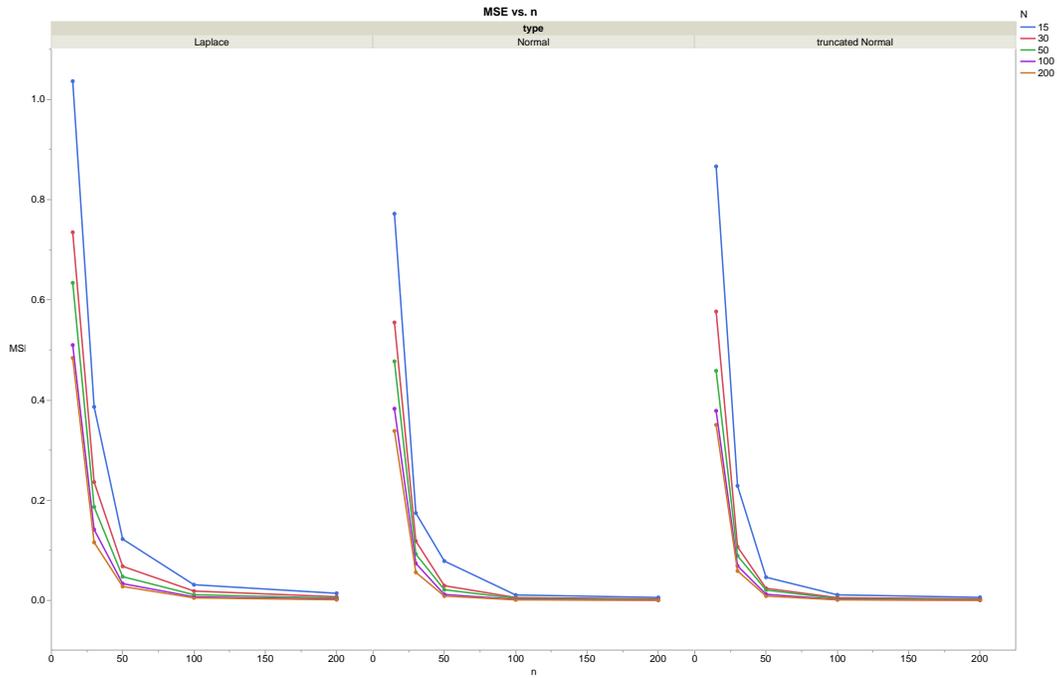}
	\caption{The MSE of STS estimates for {\it truncated Normal}, {\it Normal}, {\it Laplace} error-terms/effects with $\sigma=0.1, \lambda=0.1$.}
	\label{6.9}
\end{figure}

\begin{figure}[!h]
	\centering 
	\includegraphics[width=0.85\textwidth]{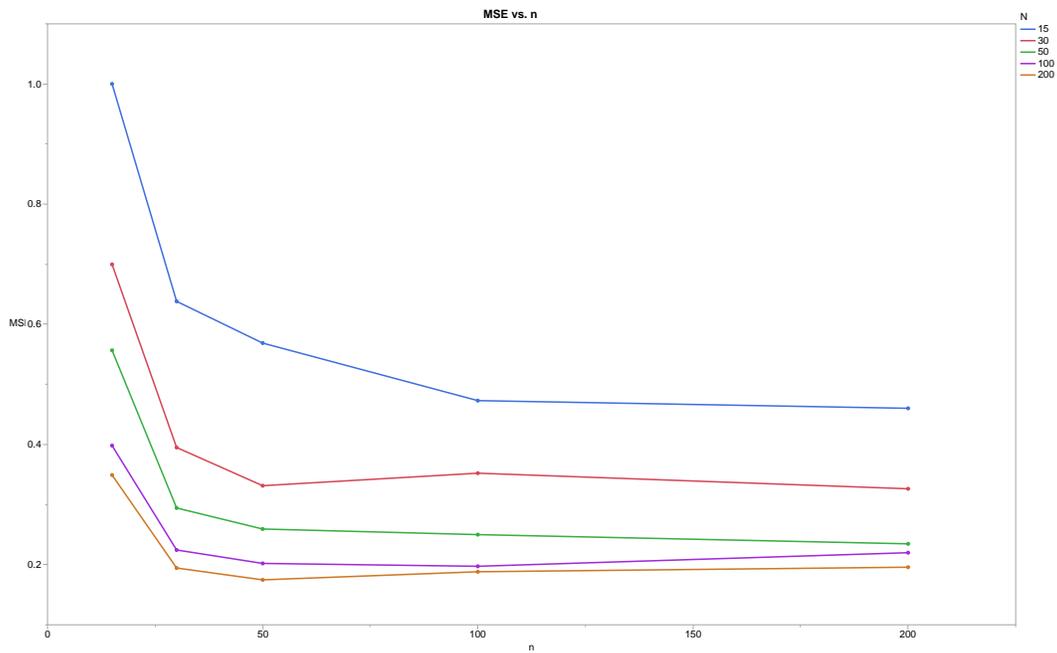}
	\caption{The MSE of STS estimates for {\it truncated Normal} error-terms/effects with $\sigma=0.05, \lambda=1$.}
	\label{6.10}
\end{figure}

\begin{figure}[!h]
	\centering 
	\includegraphics[width=0.8\textwidth]{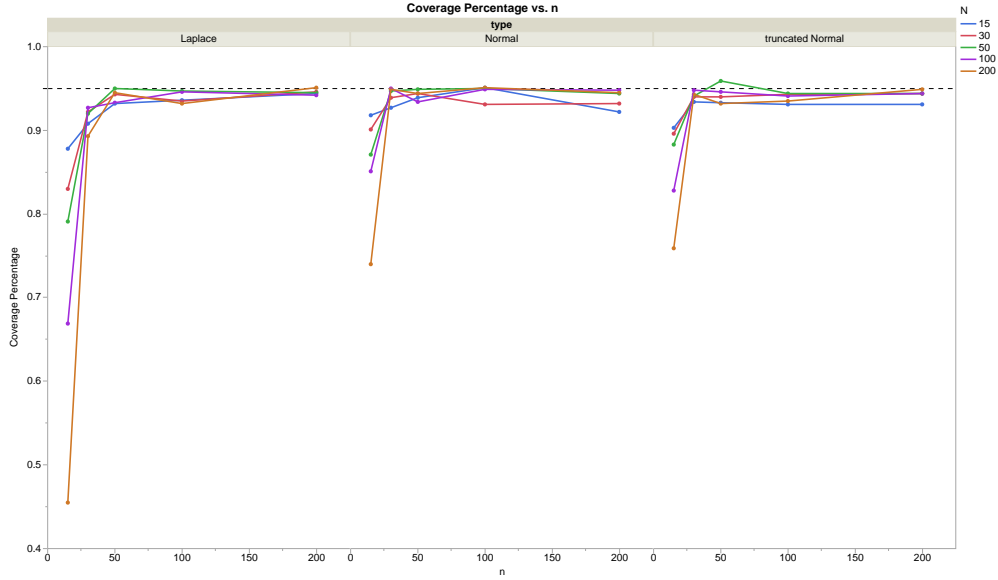}
	\caption{Coverage Percentage of the CI for the {\it truncated Normal}, {\it Normal}, {\it Laplace} error-terms/effects with $\sigma=0.5, \lambda=0.5$.}
	\label{6.13}
\end{figure}

We further considered simulations using the Laplace distributions for the error terms and random effects terms. The results are provided in Table \ref{5.16} and Table \ref{5.17}. We can see the performance of STS estimates in Laplace error terms case is consistent with normal error case. 
We also illustrate the these simulation results in Figures \ref{6.9} - \ref{6.13}. Figure \ref{6.9} depicts the MSE of STS estimates for {\it truncated Normal}, {\it Normal}, {\it Laplace} error-terms/effects with $\sigma=0.1, \lambda=0.1$. Figure \ref{6.10} depicts the MSE of STS estimates for {\it truncated Normal} error-terms/effects with $\sigma=0.05, \lambda=1$. Figure \ref{6.13} illustrate the  coverage percentage of the CI for the {\it truncated Normal}, {\it Normal}, {\it Laplace} error-terms/effects with $\sigma=0.5, \lambda=0.5$.
\begin{table}[!h]
{\absmall
\begin{center}
		\begin{tabular}{c c c c c c}
			\hline\hline
			& n=15 &n=30 &n=50 &n=100 &n=200\\		
			\hline
			N=15     &1.03613 &0.38643 &0.12267 &0.03157 &0.01450\\
			\hline
			N=30     & 0.73469 &0.23642 &0.06831 &0.01897 &0.00756\\ 
			\hline
			N=50     & 0.63382 &0.18683 &0.04771 &0.01161 &0.00492\\
			\hline
			N=100    & 0.50973 &0.14164 &0.03378 &0.00738 &0.00288\\
			\hline
			N=200   & 0.48408 &0.11612 &0.02806 &0.00532 &0.00159\\
			\hline\hline			
		\end{tabular}
		\caption{The MSE of STS estimates for {\it Laplace} error-terms/effects with $\sigma=0.1, \lambda=0.1$.}
		\label{5.16}
\end{center}
}
\end{table}
\vfill
\begin{table}[!h]
{\absmall
\begin{center}
		\begin{tabular}{c c c c c c}
			\hline\hline
			& n=15 &n=30 &n=50 &n=100 &n=200\\		
			\hline
			N=15     &  0.878 &0.908 &0.932 &0.936 &0.944 \\
			\hline
			N=30     &0.830 &0.922 &0.943 &0.935 &0.946\\ 
			\hline
			N=50     &0.791 &0.920 &0.950 &0.947 &0.945\\
			\hline
			N=100  & 0.669 &0.927 &0.933 &0.946 &0.942\\
			\hline 
			N=200   & 0.455 &0.893 &0.945 &0.932 &0.951 \\
			\hline\hline			
		\end{tabular}
		\caption{Coverage Percentage of the CI for the {\it Laplace} error-terms/effects with $\sigma=0.5, \lambda=0.5$.}
		\label{5.17}
\end{center}
}
\end{table}
\vfill

\vfill\eject

\subsection{Illustrating the {\it Recycled} STS Estimation Procedure}

In this section, we provide the results of the simulation studies corresponding to Theorem \ref{5.1.1} and \ref{5.1.2} concerning the {\it recycled}  STS estimation procedure with $\hat\theta^*_{_{STS}}$.  We considered the same compartmental model as given in the previous subsection, however again with $p=1$.  Accordingly, we choose $\theta_2$ to represent the model's unknown parameter and set,  for the simulations,  $\thzz=0.8$, for each $i$. As before, we generated the values of $\{t_{ij},j=1,\dots,n\}$ uniformly from the $[0,8]$ interval, and draw the error terms, $\e_{ij}$ and the random effects terms,  $b_i$,  from the {\it truncated Normal} distribution. 

For each simulation run, we calculated the value of  $\hat{\theta}_{_{STS}}$ as in section 4.2, then with $B=1,000$, we generated $B\times N$ independent replications of the random weights $\w_i=(w_{i1},w_{i2},\dots,w_{in})$ and $B=1,000$ independent replications of the random weight $\u=(u_1,u_2,\dots,u_N)$, to obtain $\hat{\theta}^{*1}_{_{STS}}$, $\hat{\theta}^{*2}_{_{STS}}$, $\dots$, $\hat{\theta}^{*B}_{_{STS}}$. The correspond 95\% Confidence Intervals were formed. 
With $\sigma= 1$, $\lambda= 1$ a total of $M=2000$ replications of such simulations were executed to determine the percentage of times the true value of the parameter estimates was contained in the interval and average confidence interval length was calculated. The Coverage Percentages with average confidence interval lengths are provided in Table \ref{5.18} to Table \ref{5.21}.

Table \ref{5.18} demonstrates the results of the asymptotic results of Section 4. Table \ref{5.19} to \ref{5.21} provide Coverage Percentages with average confidence interval lengths,  with random weights set to be {\it Multinomial, Dirichlet } or {\it Exponential} distributed . From these results we can see with $N$ and $n$ both increase, the Coverage Percentages converges  to 0.95 as expected. Also notice that Coverage Percentages derived from the {\it recycled} STS are more accurate (closer to 0.95) than the asymptotic result, especially when $n$ and $N$ are small.

\begin{figure}[!h]
	\centering 
	\includegraphics[width=1\textwidth]{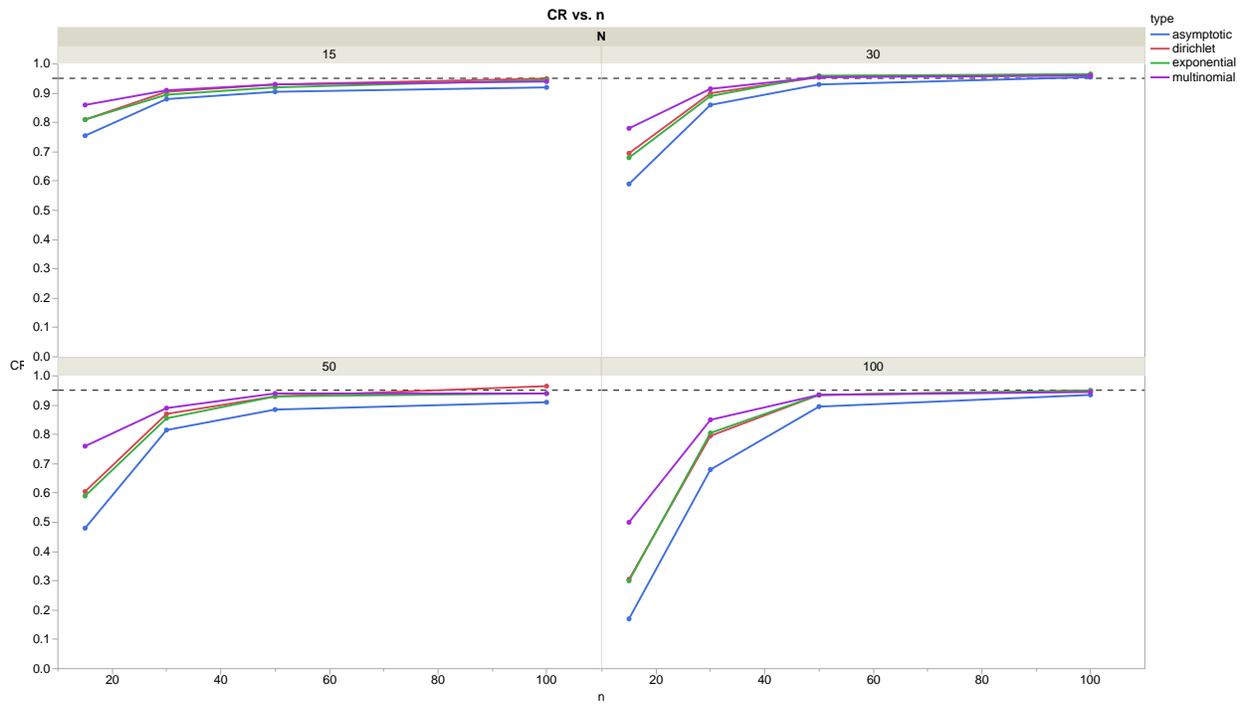}
	\caption{Coverage Percentage of the CI for the {\it truncated Normal} error-terms/effects with $\sigma=1, \lambda=1$.}
	\label{6.11}
\end{figure}

\begin{figure}[!h]
	\centering 
	\includegraphics[width=1\textwidth]{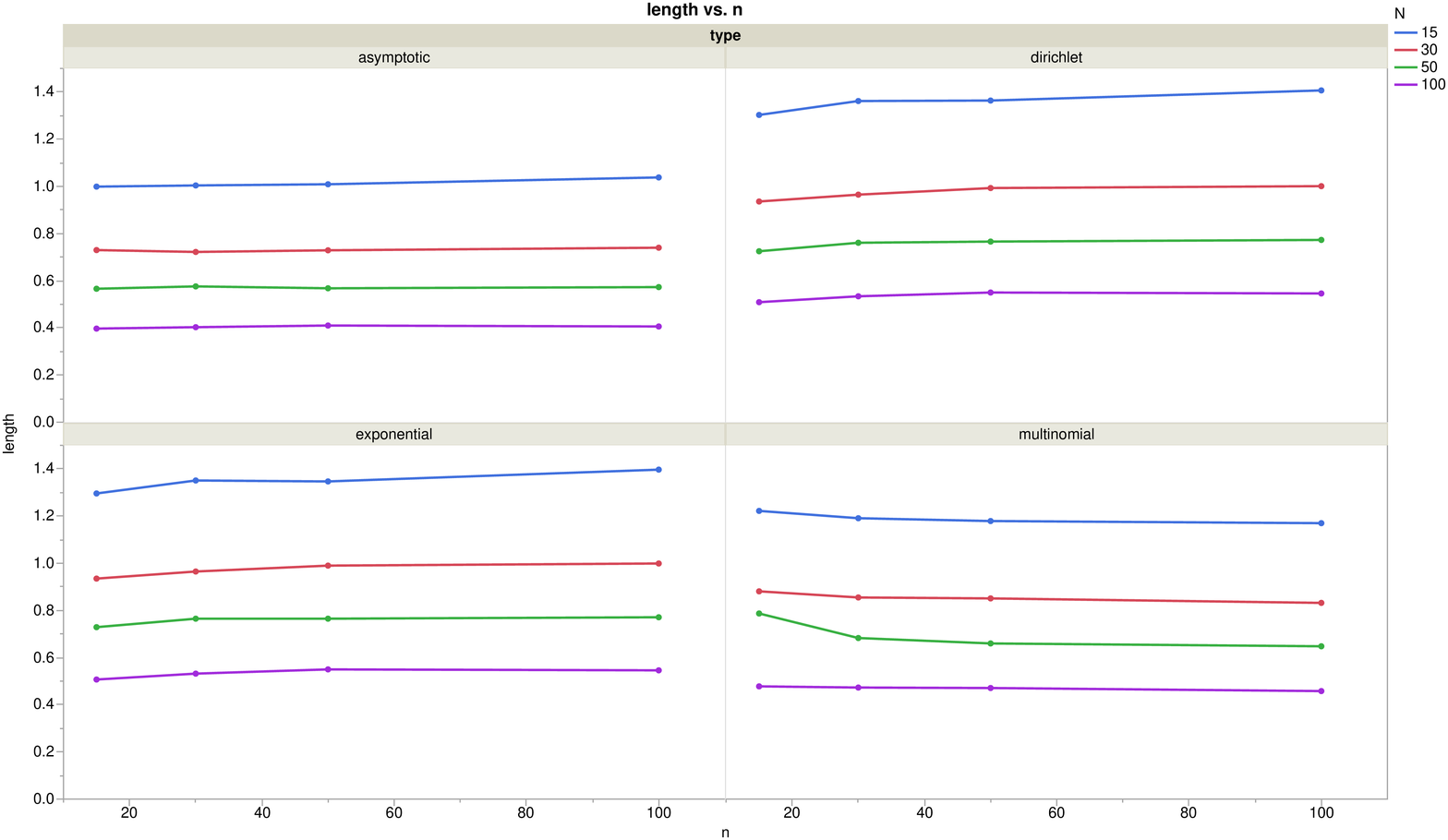}
	\caption{Average length of the CI for the {\it truncated Normal} error-terms/effects with $\sigma=1, \lambda=1$.}
	\label{6.12}
\end{figure}

\begin{table}
{\absmall
\begin{center}
		\begin{tabular}{c c c c c }
			\hline\hline
			& n=15 &n=30 &n=50 &n=100\\		
			\hline
			N=15   &0.755 &0.880 &0.905&0.920\\
			   &0.999&1.004&1.009&1.038
			\\
			\hline
			N=30     &0.590&	0.860&0.930&0.955\\
			&0.730&0.722&0.729&0.740	\\ 
			\hline
			N=50     &0.48&0.815&0.885&0.955\\
                     &0.566&0.576&0.568&0.573\\
			\hline
			N=100  & 0.170&0.680&0.895&0.935\\
			&0.397 &0.403&0.410&0.406\\
			\hline\hline			
		\end{tabular}
		\caption{Simulated Coverage Percentage of the CI for the {\it truncated Normal} error-terms/effects with $\sigma=1, \lambda=1$.}
		\label{5.18}
\end{center}
}
\end{table}
\begin{table}
{\absmall
\begin{center}
		\begin{tabular}{c c c c c }
			\hline\hline
			& n=15 &n=30 &n=50 &n=100\\		
			\hline
			N=15   &0.860&0.910&0.930&0.940\\
			   &1.222&1.191&1.179&1.170
			\\
			\hline
			N=30     &0.780&	0.915&0.955&0.960	\\ 
			         &0.881&0.855&0.851&0.832\\
			\hline
			N=50     &0.760&0.890&0.940&0.940\\
			&0.787&0.683&0.660&0.648\\
			\hline
			N=100  & 0.500&0.850&0.935&0.945\\
			&0.478&0.473&0.471&0.458\\
			\hline\hline			
		\end{tabular}
		\caption{Coverage Percentage of the CI for the {\it truncated Normal} error-terms/effects with $\sigma=1, \lambda=1$ and with {\it Multinomial} random weights.}	
		\label{5.19}
\end{center}
}
\end{table}
\begin{table}
{\absmall	
\begin{center}
		\begin{tabular}{c c c c c }
			\hline\hline
			& n=15 &n=30 &n=50 &n=100\\		
			\hline
			N=15   &0.810&0.905&0.930&0.950\\
			   &1.303&1.362&1.364&1.407\\
			\hline
			N=30     &0.695&0.900&0.955&0.965\\
			&0.936&0.965&0.993&1.001	\\ 
			\hline
			N=50     &0.605&0.870&0.930&0.965\\
			&0.725&0.761&0.766&0.773\\
			\hline
			N=100  &0.305&0.795&0.935&0.950\\
			&0.509&0.534&0.550&0.546\\
			\hline\hline			
		\end{tabular}
		\caption{Coverage Percentage of the CI for the {\it truncated Normal} error-terms/effects with $\sigma=1, \lambda=1$ and with {\it Dirichlet} random weights.}
\end{center}
}
\label{5.20}
\end{table}
\begin{table}
{\absmall
\begin{center}
		\begin{tabular}{c c c c c }
			\hline\hline
			& n=15 &n=30 &n=50 &n=100\\		
			\hline
			N=15   &0.810&0.895&0.920&0.945\\
			      &1.296&1.351&1.347&1.397\\
			\hline
			N=30     &0.680&	0.890&0.960&0.965\\
			&0.935&0.965&0.990&0.999	\\ 
			\hline
			N=50     &0.590&0.855&0.930&0.940\\
			&0.729&0.765&0.765&0.771\\
			\hline
			N=100  & 0.300&0.805&0.935&0.950\\
			&0.507&0.532&0.550&0.546\\
			\hline\hline			
		\end{tabular}
		\caption{Coverage Percentage of the CI for the {\it truncated Normal} error-terms/effects with $\sigma=1, \lambda=1$ and with {\it Exponential} random weights.}
		\label{5.21}
\end{center}
}
\end{table}
To complement of the simulations, we also considered the  Laplace distribution for the error and random effects terms and present the corresponding simulation results Tables \ref{5.22} - \ref{5.25}, below. Table \ref{5.22} demonstrates the results from asymptotic result as in Section 4. Table \ref{5.23} to \ref{5.25} present Coverage Percentages with average confidence interval lengths with weights set to be according to the {\it Multinomial, Dirichlet } and the {\it Exponential} distributions. The results have similar performance as in normal random component case. Also notice that Coverage Percentages derived from the {\it recycled}  STS method are also more accurate (closer to 0.95) than the asymptotic result, especially for smaller $n$ and $N$.. We also illustrate these simulation results in Figure \ref{6.11} and \ref{6.12}. Figure \ref{6.11} is coverage percentage of the CI for the {\it truncated Normal} error-terms/effects with $\sigma=1, \lambda=1$. Figure \ref{6.12} is average length of the CI for the {\it truncated Normal} error-terms/effects with $\sigma=1, \lambda=1$. From this figure we can observe that with an increasing $N$, the average length of the CI is decreasing, however, with only $n$ increase the length will not decrease, which is consistent with our main results.
\begin{table}
{\absmall	
	\begin{center}
			\begin{tabular}{c c c c c }
				\hline\hline
				& n=15 &n=30 &n=50 &n=100\\		
				\hline
				N=15   &0.790&0.895&0.910&0.895\\
				&0.998&0.974&0.964&1.007\\
				\hline
				N=30     &0.730&0.885&0.870&0.940\\
				&0.714&0.726&0.715&0.714
					\\ 
				\hline
				N=50     &0.475&0.840&0.925&0.940\\
				&0.559&0.562&0.546&0.552\\
				\hline
				N=100  &0.220&0.715&0.895&0.960
				\\
				&0.395&0.388&0.390&0.397
				\\
				\hline\hline			
			\end{tabular}
		\caption{Simulated Coverage Percentage of the CI for the {\it Laplace} error-terms/effects with $\sigma=1, \lambda=1$.}	
		\label{5.22}
	\end{center}
	}
\end{table}
  
\begin{table}
	{\absmall
	\begin{center}
			\begin{tabular}{c c c c c }
				\hline\hline
				& n=15 &n=30 &n=50 &n=100\\		
				\hline
				N=15   &0.885&0.935&0.950&0.950
				\\
				&1.205&1.182&1.160&1.171
				
				\\
				\hline
				N=30     &0.905&0.960&0.915&0.955\\ 
				&0.865&0.854&0.846&0.815
				\\
				\hline
				N=50     &0.760&0.930&0.965&0.960
				\\
				&0.677&0.670&0.653&0.637
				\\
				\hline
				N=100  &0.620&0.825&0.935&0.965
				\\
				&0.475&0.465&0.459&0.456
				\\
				\hline\hline			
			\end{tabular}
				\caption{Coverage Percentage of the CI for the {\it Laplace} error-terms/effects with $\sigma=1, \lambda=1$ and with {\it Multinomial} random weights.}
				\label{5.23}
	\end{center}
	}
\end{table}
\begin{table}
	{\absmall	
	\begin{center}
			\begin{tabular}{c c c c c }
				\hline\hline
				& n=15 &n=30 &n=50 &n=100\\		
				\hline
				N=15   &0.830&0.930&0.960&0.965
				\\
				&1.309&1.350&1.367&1.422
				\\
				\hline
				N=30     &0.815&0.935&0.910&0.965
				\\
				&0.926&0.974&0.984&0.980
					\\ 
				\hline
				N=50     &0.615&0.915&0.965&0.965
				\\
				&0.721&0.758&0.757&0.768
				\\
				\hline
				N=100  &0.440&0.800&0.940&0.985
				\\
				&0.508&0.528&0.537&0.546
				\\
				\hline\hline			
			\end{tabular}
			\caption{Coverage Percentage of the CI for the {\it Laplace} error-terms/effects with $\sigma=1, \lambda=1$ and with {\it Dirichlet} random weights.}
		\label{5.24}	
	\end{center}
	}
\end{table}
\begin{table}
	{\absmall
	\begin{center}
			\begin{tabular}{c c c c c }
				\hline\hline
				& n=15 &n=30 &n=50 &n=100\\		
				\hline
				N=15   &0.845&0.930&0.950&0.960
				\\
				&1.302&1.334&1.355&1.407
				\\
				\hline
				N=30     &0.820&0.940&0.935&0.965
				\\
				&	0.923&0.969&0.982&0.979
				\\ 
				\hline
				N=50     &0.600&0.910&0.965&0.955
				\\
				&0.717&0.757&0.757&0.764
				\\
				\hline
				N=100  &0.435&0.815&0.945&0.985
				 \\
				&0.507&0.526&0.537&0.544
				\\
				\hline\hline			
			\end{tabular}
			\caption{Coverage Percentage of the CI for the {\it Laplace} error-terms/effects with $\sigma=1, \lambda=1$ and with {\it Exponential} random weights.}	
			\label{5.25}
	\end{center}
	}
\end{table}
\vfill

\section{Technical Details and Proofs}

\subsection{Technical Details and Proofs -- the STS Estimation Case}\label{5.33}
In this section we provide the technical results needed for the proofs of Theorems \ref{4.3.1} and \ref{4.3.2} on the STS estimator $\hat\theta_{_{STS}}$  in the hierarchical nonlinear regression model.  In the sequel, we let $\phi_{1ij}(\thtt):=\phi_{ij}^{'}(\thtt)$ (see (\ref{4.1.1.1})), and set $K$ to denote a {\it generic} constant.  Recall that (see {\it Assumption A(1)}), 
$$
a_{n_i}^{2}:=\sigma^2\jsum E(f^{'2}_{ij}(\thzz+b_i))\to \infty\ \ as \ \ n_i\to \infty.
$$

\begin{lemma}\label{7.3.1}
	Under the conditions of {\it Assumption A}, for some  $K>0$
	\[
	a_{n_i}^{-2}\underset{|t|\le K}{\sup}\jsum \phi_{1ij}(b_{i1})-\frac{1}{\sigma^2}\to 0 \ \ a.s.,
	\]
	where $b_{i1}:=b_{1n_{i}}(t)$ is a sequence such that $\underset{|t|\le K}{\sup}|b_{i1}-b_i-\thzz|\to 0,\ \ \ a.s.,$
as $n_i\to\infty$.
\end{lemma}

{\it \underbar{Proof of Lemma \ref{7.3.1}}:} Since $\phi_{1ij}(\thtt):=\phi_{ij}^{'}(\thtt)$, we have 
\[
\phi_{1ij}(\thtt)\equiv f^{'2}_{ij}(\thtt)-\e_{ij}f_{ij}^{''}(\thtt)-(f_{ij}(\thzz+b_i)-f_{ij}(\thtt))f_{ij}^{''}(\thtt).
\]
Accordingly,  we first note that,
\beqn
\left|a_{n_i}^{-2}\underset{|t|\le K}{\sup}\jsum \phi_{1ij}(b_{i1})-\frac{1}{\sigma^2}\right|&\le&\left|a_{n_i}^{-2}\underset{|t|\le K}{\sup}\jsum f^{'2}_{ij}(b_{i1})-\frac{1}{\sigma^2}\right|\\
&+&a_{n_i}^{-2}\underset{|t|\le K}{\sup}\left|\jsum\e_{ij}f_{ij}^{''}(b_{i1})\right|\\
&+&a_{n_i}^{-2}\underset{|t|\le K}{\sup}\left|\jsum(f_{ij}(\thzz+b_i)-f_{ij}(b_{i1}))f_{ij}^{''}(b_{i1})\right|.
\eeqn
By {\it Assumption A}  $(3)$, we have 
$
a_{n_i}^{-2}\underset{|t|\le K}{\sup}\jsum f^{'2}_{ij}(b_{i1})-\frac{1}{\sigma^2}\to 0 \ \ a.s.,
$
and by {\it Assumption A $(2)$}   and Corollary A in Wu (1981), we also have, 
\[
a_{n_i}^{-2}\underset{|t|\le K}{\sup}\left|\jsum\e_{ij}f_{ij}^{''}(b_{i1})\right|\to 0\ \ a.s..
\]
Finally, the last term converge to $0$ a.s. by {\it Assumption A}, an application of Cauchy-Schwarz inequality and Corollary A in Wu (1981). 
Thus we have 
\[
a_{n_i}^{-2}\underset{|t|\le K}{\sup}\jsum \phi_{1ij}(b_{i1})-\frac{1}{\sigma^2}\to 0 \ \ a.s..
\].
\qed
\begin{lemma}\label{7.3.2}
	Let $X_i$ be a sequence of random variables bounded in probability
	and let $Y_i$ be a sequence of random variables which satisfies $\frac{1}{n}\sum_{i=1}^{n} |Y_i|\to 0$ in probability.
	Then	
	$
	\frac{1}{n}\sum_{i=1}^{n} X_iY_i\overset{p}{\to}0.
$
\end{lemma}
{\it \underbar{Proof of Lemma \ref{7.3.2}}:}
Since $X_i$ is bounded in probability, for any $\e>0$, there is $K_\e$ such that with sufficient large i,
$
P(|X_i|>K_\e)<\e.
$
Then
\beqn
\underset{n\to\infty}{\lim}P(|\frac{1}{n}\sum_{i=1}^{n} X_iY_i|>\e)&=&\underset{n\to\infty}{\lim}\l P(|\frac{1}{n}\sum_{i=1}^{n} X_iY_i|>\e, |X_i|<K_\e)\r\\
&+&\underset{n\to\infty}{\lim}\l P(|\frac{1}{n}\sum_{i=1}^{n} X_iY_i|>\e, |X_i|>K_\e)\r\\
&\le&\underset{n\to\infty}{\lim}P(\frac{1}{n}\sum_{i=1}^{n} |\frac{X_i}{K_\e}Y_i|>\frac{\e}{K_\e}, |X_i|<K_\e)+\e\\
&\le&\underset{n\to\infty}{\lim}P(\frac{1}{n}\sum_{i=1}^{n} |Y_i|>\frac{\e}{K_\e}, |X_i|<K_\e)+\e=\e, 
\eeqn
from which the desired result follows.
\qed
\begin{lemma}\label{7.3.4}
There exists a $K>0$ such that for any $\e>0$, for any $i$,
\[
P\l\left|a_{n_i}^{-1}\jsum \phi_{ij}(\thzz+b_i)\right|>K\r<\frac{\e}{2}.
\]	
\end{lemma}
{\it \underbar{Proof of Lemma \ref{7.3.4}}:} Since $\e_{ij}$ and $b_i$ are independent, for each $i=1, \dots, N$, we have that for any $j_1\ne j_2$, 
\beqn
E(\phi_{ij_1}(\thzz+b_i)\phi_{ij_2}(\thzz+b_i))&=&E[E(\phi_{ij_1}(\thzz+b_i)\phi_{ij_2}(\thzz+b_i)|b_i)]\\
&=&E[E(\e_{ij_1}\e_{ij_2}f^{'}_{ij_1}(\thzz+b_i)f^{'}_{ij_2}(\thzz+b_i)|b_i)]\\
&=&E[E(\e_{ij_1})E(\e_{ij_2})f^{'}_{ij_1}(\thzz+b_i)f^{'}_{ij_2}(\thzz+b_i)]\\
&=&0.
\eeqn
Similarly, 
\beqn
E(\phi_{ij_1}(\thzz+b_i))&=&E[E(\e_{ij_1}f^{'}_{ij_1}(\thzz+b_i)|b_i)]\\
&=&E[E(\e_{ij_1})f^{'}_{ij_1}(\thzz+b_i)]\\
&=&0.
\eeqn    
Hence, we have, 
\[
E(\phi_{ij_1}(\thzz+b_i)\phi_{ij_2}(\thzz+b_i))=E(\phi_{ij_1}(\thzz+b_i))E(\phi_{ij_2}(\thzz+b_i)).
\]
To conclude that, 
\beqn
Var\ll\jsum \phi_{ij}(\thzz+b_i)\rr&=&\jsum Var(\phi_{ij}(\thzz+b_i))\\
&=&\jsum Var(\e_{ij}f^{'}_{ij}(\thzz+b_i))\\
&=&\jsum E(\e^2_{ij})E(f^{'2}_{ij}(\thzz+b_i))\\
&=&\sigma^2\jsum E(f^{'2}_{ij}(\thzz+b_i))\equiv a_{n_i}^2.
\eeqn
Accordingly, there exists a $K>0$ such that for any $\e>0$, for any $i$,
\[
P\l\left|a_{n_i}^{-1}\jsum \phi_{ij}(\thzz+b_i)\right|>K\r<\frac{\e}{2}.
\]

\qed

\noindent{\underline{\it Proof of Theorem \ref{4.3.1}:}} 
Let 
\be\label{7.9}
S_{n_i}(t):= a_{n_i}^{-1}\jsum \l\phi_{ij}(\thzz+b_i+a_{n_i}^{-1}t)-\phi_{ij}(\thzz+b_i)\r-\frac{t}{\sigma^2}.
\ee
Next we will show for any given constant $K$,
\be\label{3.1}
\underset{|t|\le K}{sup}|S_{n_i}(t)|\to 0\ \ a.s.
\ee
By a Taylor expansion,
$
\phi_{ij}(\thzz+b_i+a_{n_i}^{-1}t)=\phi_{ij}(\thzz+b_i)+\phi_{1ij}(b_{i1})a_{n_i}^{-1}t, \ \ 
$
where $b_{i1}=\thzz+b_i+ca_{n_i}^{-1}t$ for some $0<c<1$. Accordingly we obtain that, 
\beqn
\underset{|t|\le K}{\sup}|S_{n_i}(t)|&=&\underset{|t|\le K}{\sup} \left|a_{n_i}^{-1}\jsum \phi_{1ij}(b_{i1})a_{n_i}^{-1}t-\frac{t}{\sigma^2} \right|\\
&=&K\left| a_{n_i}^{-2}\underset{|t|\le K}{\sup}\jsum\phi_{1ij}(b_{i1})-\frac{1}{\sigma^2}\right|.\\
\eeqn
By Lemma \ref{7.3.1},
$
a_{n_i}^{-2}\underset{|t|\le K}{\sup}\jsum\phi_{1ij}(b_{i1})-\frac{1}{\sigma^2}\to0\ \ a.s.
$
Thus,  we have proved (\ref{3.1}). Next, by (\ref{7.9}),
\beqn
A_{n_i}(t):=a_{n_i}^{-1}t\jsum \phi_{ij}(\thzz+b_i+a_{n_i}^{-1}t)= tS_{n_i}(t)+a_{n_i}^{-1}t\jsum \phi_{ij}(\thzz+b_i)+\frac{t^2}{\sigma^2}.
\eeqn
Thus,
\beqn
\underset{|t|=K}{\inf}A_{n_i}(t)\ge -K\underset{|t|=K}{\sup}|S_{n_i}(t)|-Ka_{n_i}^{-1}\left|\jsum \phi_{ij}(\thzz+b_i)\right|+\frac{K^2}{\sigma^2}.
\eeqn
By lemma \ref{7.3.4} there exists a $K>0$ such that for any $\e>0$, for any $i$,
\be\label{3.2}
P\l\left|a_{n_i}^{-1}\jsum \phi_{ij}(\thzz+b_i)\right|>K\r<\frac{\e}{2}.
\ee
So that by (\ref{3.2}) and (\ref{3.1}) we may choose $K$ large enough such that for sufficiently large $n_i$,
\beqn
P(\underset{|t|=K}{\inf}A_{n_i}(t)\ge0)&\ge& P(\underset{|t|=K}{\sup}|S_{n_i}(t)|+a_{n_i}^{-1}\left|\jsum \phi_{ij}(\thzz+b_i)\right|\le \frac{K}{\sigma^2})\\
&=&1-P(\underset{|t|=K}{\sup}|S_{n_i}(t)|+a_{n_i}^{-1}\left|\jsum \phi_{ij}(\thzz+b_i)\right|>\frac{K}{\sigma^2})\\
&\ge&1-P(\underset{|t|=K}{\sup}|S_{n_i}(t)|>\frac{K}{4\sigma^2})-P(a_{n_i}^{-1}\left|\jsum \phi_{ij}(\thzz+b_i)\right|>\frac{K}{4\sigma^2})\\
&\ge& 1-\e.
\eeqn
By the continuity of $\jsum \phi_{ij}(\thtt)$ in $\thtt$, we have,  for sufficiently large $n_i$, that there exists a constant $K$ such that the equation 
\[
\jsum \phi_{ij}(\thzz+b_i+a_{n_i}^{-1}t)=0, 
\]
has a root $t=T_{ni}$ in $|t|\le K$ with probability larger than $1-\e$. That is,  we have
\[
\thnn=\thzz+b_i+a_{ni}^{-1}T_{ni}, 
\]
where $|T_{ni}|<K$ in probability. Thus, by Lemma \ref{7.3.2},
\[
\hat\theta_{_{STS}}-\thzz=\frac{1}{N}\isum b_i +\frac{1}{N}\isum a_{ni}^{-1}T_{ni}\overset{p}{\to}0.
\]
\qed

\noindent For establishing the asymptotic normality result as stated in Theorem \ref{4.3.2}, we need the following Lemma.

\begin{lemma}\label{7.3.3}
	Under the conditions of {\it Assumptions A},
	\[
	\frac{1}{\sqrt{N}}\isum a_{n_i}^{-2}\jsum\phi_{ij}(\thzz+b_i)\overset{p}{\to} 0.
	\]
\end{lemma}
{\it \underbar{Proof of Lemma \ref{7.3.3}}:} Let
$
X_{ni}:= a_{n_i}^{-1}\jsum \phi_{ij}(\thzz+b_i), 
$
where,  by proof of Theorem \ref{4.3.1} we have 
$E(X_{ni})=0$ and $Var(X_{ni})=1$.  Thus,
\[
\frac{1}{\sqrt{N}}\isum a_{n_i}^{-2}\jsum\phi_{ij}(\thzz+b_i)= \frac{1}{\sqrt{N}}\isum a_{n_i}^{-1}X_{ni}.
\]
Now,  for any $\e>0$,
\beqn
P(\left|\frac{1}{\sqrt{N}}\isum a_{n_i}^{-1}X_{ni}\right|>\e)\le \frac{\isum \frac{1}{a_{ni}^2}}{N\e^2}\to 0.
\eeqn
Accordingly,  we have
$
\frac{1}{\sqrt{N}}\isum a_{n_i}^{-2}\jsum\phi_{ij}(\thzz+b_i)\overset{p}{\to} 0, 
$
as required. 
\qed
\bigskip

\noindent{\underline{\it Proof of Theorem \ref{4.3.2}:} } We first note that by Lemma \ref{7.3.1} and ($\ref{7.9}$),
\[
\thnn-\thzz-b_i=-a_{ni}^{-2}\sigma^2\jsum\phi_{ij}(\thzz+b_i)- a_{n_i}^{-1}\sigma^2S_{n_i}(T_{ni}). 
\]
Thus,
\[
\hat\theta_{_{STS}}-\thzz= \frac{1}{N}\isum b_i -\frac{\sigma^2}{N}\isum a_{n_i}^{-2}\jsum\phi_{ij}(\thzz+b_i)-\frac{\sigma^2}{N}\isum a_{n_i}^{-1}S_{n_i}(T_{ni}). 
\]
Recall that $\isum b_i/N\to E(b_1)\equiv 0$. In view of (\ref{3.1}) and since, $\underset{N,ni\to \infty}{\lim} N/a_{ni}^2<\infty$,   we have 
\[
\frac{\sigma^2}{\sqrt{N}}\isum a_{n_i}^{-1}S_{n_i}(T_{ni})\to 0 \ \ a.s..
\]
Finally, from Lemma \ref{7.3.3},
$$
\frac{1}{\sqrt{N}}\isum a_{n_i}^{-2}\jsum\phi_{ij}(\thzz+b_i)\overset{p}{\to} 0.
$$
Thus, it follows that $
\lambda^{-1}\sqrt{N}(\hat\theta_{_{STS}}-\thzz)\Rightarrow {\cal {N}}(0,1).
$
\qed

\subsection{Technical Details and Proofs -- the Recycled STS Estimation Case}\label{5.44}

In this section we provide the technical results needed for the proofs of Theorems \ref{5.1.1} and \ref{5.1.2} on the {\it recycled} STS estimator,  $\hat\theta^*_{_{STS}}$,  in the hierarchical nonlinear regression model. We begin with a re-statement of Lemma 2 from Boukai and Zhang (2018) which is concerned with the general random weights under {\it Assumption W}. . 

\begin{lemma}\label{2.4.3} Let $\w_n=(w_{1:n}, w_{1:n},  \dots, w_{n:n})^\t$ be random weights that satisfy the conditions of {\it Assumption W}. Then With $W_{i}=(w_{i:n}-1)/\tau_n, \ i=1\dots, n$ and $\bar W_n:=\frac{1}{n}\ssum W_i$ we have,  as $n\to \infty$, that $(i)\ \ \frac{1}{n}\ssum W_i\overset{p^*}{\to}0$  $(ii)\ \ \frac{1}{n}\ssum W_i^2\overset{p^*}{\to}1$
	and hence $(iii)\ \ \frac{1}{n}\ssum (W_i-\bar{W}_n)^2\overset{p^*}{\to}1$. 
\end{lemma}

\begin{lemma}\label{2.4.6}
		Under the conditions of {\it Assumption W},
$
	\frac{1}{n}\ssum w_{i:n}-1\overset{p^*}{\to}0 ,
$
Further, let $\U_n=(u_1,u_2, \dots, u_n)^\t$ denote a vector of $n$ $i.i.d$ random variables that is independent of $\w_n$ with $E(u_i)=0$,  $E(u_i^2)<\infty$. Then, conditional on the given value of the $\U_n$, we have $\frac{1}{n}\ssum u_i w_{i:n} \overset{p^*}{\to}0$, as $n\to \infty$. 

\end{lemma}

\noindent {\it \underbar{Proof of Lemma \ref{2.4.6}}:} We first note that 
\beqn
	E^*(\frac{1}{n}\ssum (w_{i:n}-1))^2&=&E^*(\frac{\tau_n}{n}\ssum W_i)^2\\
&=&\frac{\tau_n^2}{n^2}\ssum E^*(W_i^2)+\frac{\tau_n^2}{n^2}\underset{i_1\ne i_2}{\sum}E^*(W_{i_1}W_{i_2})\\
&=&\frac{\tau_n^2}{n}+\frac{\tau_n^2}{n^2}n(n-1)O(\frac{1}{n})\to 0, \ \ \ \text{as} \ \ n\to \infty.
\eeqn
To conclude that, $\frac{1}{n}\ssum w_i-1\overset{p^*}{\to}0$, as $n\to \infty$. As for the second assertion, we note that since 
\[
\frac{1}{n}\ssum u_iw_{i:n}=\frac{\tau_n}{n}\ssum u_iW_i+\frac{1}{n}\ssum u_i, 
\]
and since $\ssum u_i/n\to 0$, as $n\to\infty$, we may only consider the first term. To that end, we note that
\beqn
E^*(\frac{\tau_n}{n}\ssum u_iW_i)^2&=&\frac{\tau_n^2}{n^2}\ssum E^*(u_i^2W_i^2)+\frac{\tau_n^2}{n^2}\underset{i_1\ne i_2}{\sum}E^*(W_{i_1}W_{i_2}u_{i_1}u_{i_2})\\
&\le&\l1+(n-1)O(\frac{1}{n})\r\frac{\tau_n^2}{n^2}\ssum u_i^2 \to 0,
\eeqn
as $n\to \infty$. We therefore conclude that $\frac{1}{n}\ssum u_iw_{i:n}\overset{p^*}{\to}0, $ as required. \qed
\begin{lemma}\label{7.4.1}
	Under the conditions of {\it Assumptions A} and {\it  B}, we have  that  $a_{n_i}^{-2}\jsum\phi^2_{ij}(\thnn)\overset{p}{\to} 1, $ 
for all $i=1,2,\dots,N$. 
\end{lemma}
{\it \underbar{Proof of Lemma \ref{7.4.1}}:} Since $\thnn\overset{p}{\to}\thzz$, we have
\beqn
a_{n_i}^{-2}\jsum\phi^2_{ij}(\thnn)&=&a_{n_i}^{-2}\jsum(y_{ij}-f_{ij}(\thnn))^2f^{'2}_{ij}(\thnn)\\
&=&a_{n_i}^{-2}\jsum \e_{ij}^2f_{ij}^{'2}(\thnn)+a_{n_i}^{-2}\jsum (f_{ij}(\thzz+b_i)-f_{ij}(\thnn))^2f^{'2}_{ij}(\thnn)\\
&+&2a_{n_i}^{-2}\jsum \e_{ij}(f_{ij}(\thzz+b_i)-f_{ij}(\thnn))f^{'2}_{ij}(\thnn)\\
& \equiv &B_1 +B_2 +B_3.
\eeqn
Write, 
\beqn
B_1=a_{n_i}^{-2}\jsum (\e_{ij}^2-\sigma^2)f^{'2}_{ij}(\thnn)+a_{n_i}^{-2}\sigma^2\jsum f^{'2}_{ij}(\thnn).
\eeqn
The first term in $B_1$ converges to $0$ by {\it Assumption A (3)}, and Corollary A of Wu (1981) while the second term in $B_1$ converges to 1 by {\it Assumption A (3)}. Hence $B_1\overset{p}{\to} 1$.  As for the second and third terms,  $B_2$ and $B_3$, it follows by a direct application of the Cauchy-Schwarz inequality ogether with {\it Assumption B (1)}, that  $B_2\overset{p}{\to} 0$ and $B_3\overset{p}{\to} 0$.  Accordingly, it follows that $
a_{n_i}^{-2}\jsum\phi^2_{ij}(\thnn)\overset{p}{\to} 1, 
$
as required. 
\qed
\begin{lemma}\label{7.4.2}
	Under the conditions of {\it Assumptions A} and {\it B}, for all $i$,
	\[
	E^*\big[\tau_{n_i}a_{n_i}^{-2}\underset{|t|\le K\tau_{n_i}}{\sup}\jsum W_{ij} \phi_{1ij}(b^*_{i1})\big]^2\to 0
	\]
	where $b^*_{i1}=\thnn+ca_{n_i}^{-1}t$ for some $0<c<1$, as $n_i\to\infty$.
\end{lemma}
{\it \underbar{Proof of Lemma \ref{7.4.2}}:} 
We first note that since by Theorem $\ref{4.3.1}$, we have $\thnn-b_i-\thzz\overset{p}{\to}0$,  and since
\beqn
|b^*_{i1}-b_i-\thzz|&=&|\thnn-b_i-\thzz+ca_{n_i}^{-1}t|\\
&\le&|\thnn-b_i-\thzz|+\frac{c\tau_{n_i}}{\sqrt{n_i}}\frac{\sqrt{n_i}}{a_{n_i}}\frac{|t|}{\tau_{n_i}},
\eeqn 
it follows under {\it Assumption B (3)}  that  with  $|t|\le K\tau_{n_i}$,  we have
$
b^*_{i1}-b_i-\thzz\overset{p}{\to}0.
$ 
Thus,
\beqn
&&E^*[\tau_{n_i}a_{n_i}^{-2}\underset{|t|\le K\tau_{n_i}}{\sup} \jsum W_{ij} \phi_{1ij}(b^*_{i1})]^2\\
&\le&\tau_{n_i}^2a_{n_i}^{-4}\underset{|t|\le K\tau_{n_i}}{\sup}[\jsum \phi^2_{1ij}(b^*_{i1})+O(\frac{1}{n_i})\sum_{j_1\ne j_2}\phi_{1ij_1}(b^*_{i1})\phi_{1ij_2}(b^*_{i1})]\\
&\le&\tau_{n_i}^2a_{n_i}^{-4}\underset{|t|\le K\tau_{n_i}}{\sup}[\jsum \phi^2_{1ij}(b^*_{i1})+O(\frac{1}{n_i})(n_i-1)\jsum \phi^2_{1ij}(b^*_{i1})]\\
&=&\tau_{n_i}^2a_{n_i}^{-4}[O(\frac{1}{n_i})(n_i-1)+1]\underset{|t|\le K\tau_{n_i}}{\sup}\jsum \phi^2_{1ij}(b^*_{i1}).
\eeqn
In light  of {\it Assumption B (2-3)} , and that $\tau_{n_i}^2/n_i\to0$, we only need to show, in order to complete the prrof of  Lemma \ref{7.4.2}, that   
\[
\underset{n_i\to \infty}{\lim} \ \ a_{n_i}^{-2}\underset{|t|\le K\tau_{n_i}}{\sup}\jsum \phi^2_{1ij}(b^*_{i1})<\infty.
\]
Toward that end, we note that, 
\beqn
&& a_{n_i}^{-2}\underset{|t|\le K\tau_{n_i}}{\sup}\jsum \phi^2_{1ij}(b^*_{i1})\\
&=&a_{n_i}^{-2}\underset{|t|\le K\tau_{n_i}}{\sup}\jsum(f^{'2}_{ij}(b^*_{i1})-(y_{ij}-f_{ij}(b^*_{i1}))f^{''}_{ij}(b^*_{i1}))^2\\
&\le&a_{n_i}^{-2}\underset{|t|\le K\tau_{n_i}}{\sup}\jsum f^{'4}_{ij}(b^*_{i1})+a_{n_i}^{-2}\underset{|t|\le K\tau_{n_i}}{\sup}\jsum(y_{ij}-f_{ij}(b^*_{i1}))^2f^{''2}_{ij}(b^*_{i1})\\
&+&2a_{n_i}^{-2}\underset{|t|\le K\tau_{n_i}}{\sup}\left|\jsum f^{'2}_{ij}(b^*_{i1})(y_{ij}-f_{ij}(b^*_{i1}))f^{''}_{ij}(b^*_{i1})\right|\\
& \equiv &I_{1} +I_2 +I_3.
\eeqn
It is straight forward to see that by {\it Assumption B (1)}, $\underset{n_i\to \infty}{\lim}I_1<\infty$, and that by Cauchy-Schwarz inequality $\underset{n_i\to \infty}{\lim}I_3<\infty$. Finally we write 
\beqn
I_2&=&a_{n_i}^{-2}\underset{|t|\le K\tau_{n_i}}{\sup}\jsum(\e_{ij}^2-\sigma^2)f^{''2}_{ij}(b^*_{i1})+a_{n_i}^{-2}\underset{|t|\le K\tau_{n_i}}{\sup}\jsum\sigma^2f^{''2}_{ij}(b^*_{i1})\\
&+&a_{n_i}^{-2}\underset{|t|\le K\tau_{n_i}}{\sup}\jsum(f_{ij}(\thzz+b_i)-f_{ij}(b^*_{i1}))^2f^{''2}_{ij}(b^*_{i1})\\
&+&2a_{n_i}^{-2}\underset{|t|\le K\tau_{n_i}}{\sup}\left|\jsum\e_{ij}(f_{ij}(\thzz+b_i)-f_{ij}(b^*_{i1}))f^{''2}_{ij}(b^*_{i1})\right|.
\eeqn
The first term converges to 0 in probability by {\it Assumption B (2)} and Corollary A of Wu (1981). Then, according to {\it  Assumption A (2)},
\[
\underset{n_i\to \infty}{\lim} a_{n_i}^{-2}\underset{|t|\le K\tau_{n_i}}{\sup}\jsum\sigma^2f^{''2}_{ij}(b^*_{i1}) <\infty.
\]
The third term in  $I_2$  converges to $0$ in probability by an application of the Cauchy-Schwarz inequality combined with {\it Assumption B (1) \& (2)}. Finally,  the fourth term in $I_2$, converges to $0$ in probability again, by an application  of the Cauchy-Schwarz inequality. Thus we have $\underset{n_i\to \infty}{\lim}I_2<\infty$. Accordingly, we have established that as $n_i\to\infty$,
\[
E^*\left[\tau_{n_i}a_{n_i}^{-2}\underset{|t|\le K\tau_{n_i}}{\sup}\jsum W_{ij} \phi_{1ij}(b^*_{i1})\right]^2\to 0.
\]
\qed
\begin{lemma}\label{7.4.3}
	Under the conditions of {\it Assumptions A} and {\it B}, there exists a $K>0$ such that for any $\e>0$, 
\[
P^*\l\left|a_{n_i}^{-1}\jsum W_{ij}\phi_{ij}(\thnn)\right|>K\r<\frac{\e}{2}.
\]
\end{lemma}
{\it \underbar{Proof of Lemma \ref{7.4.3}}:} 
By Lemma \ref{7.4.1},
\beqn
&& V^*(a_{ni}^{-1}\jsum W_{ij}\phi_{ij}(\thnn))\\
&=&a_{ni}^{-2}\jsum \phi^2_{ij}(\thnn)+a_{ni}^{-2}O(\frac{1}{n_i})\underset{j_1\ne j_2}{\sum} \phi_{ij_1}(\thnn)\phi_{ij_2}(\thnn)\\
&=&a_{ni}^{-2}\jsum \phi^2_{ij}(\thnn)+a_{ni}^{-2}O(\frac{1}{n_i})(\jsum \phi_{ij}(\thnn))^2
-a_{ni}^{-2}O(\frac{1}{n_i})\jsum \phi^2_{ij}(\thnn)\\
&\le&a_{ni}^{-2}(1-O(\frac{1}{n_i}))\jsum \phi^2_{ij}(\thnn)\overset{p}{\to} 1 .
\eeqn
Hence we obtain, 
\[
P^*(\left|a_{ni}^{-1}\jsum W_{ij}\phi_{ij}(\thnn)\right|>K)\le \frac{V^*(a_{ni}^{-1}\jsum W_{ij}\phi_{ij}(\thnn))}{K^2}\overset{p}{\to}\frac{1}{K^2}.
\]
Accordingly, there exists a $K>0$ such that for any $\e>0$, 
\[
P^*\l\left|a_{n_i}^{-1}\jsum W_{ij}\phi_{ij}(\thnn)\right|>K\r<\frac{\e}{2}.
\]
\qed

\noindent{\underline{\it Proof of Theorem \ref{5.1.1}:}} Let 
\be\label{7.12}
S^*_{n_i}(t):= a_{n_i}^{-1}\jsum w_{ij}\l\phi_{ij}(\thnn+a_{n_i}^{-1}t)-\phi_{ij}(\thnn)\r-\frac{t}{\sigma^2}.
\ee
First,  we will  show that for any given $K>0$,
\be\label{3.3}
E^*\left[\tau_{n_i}^{-1}\underset{|t|\le K\tau_{n_i}}{\sup}|S^*_{n_i}(t)|\right]^2\overset{p^*}{\to}0.
\ee
By a Taylor expansion we have that $
\phi_{ij}(\thnn+a_{n_i}^{-1}t)=\phi_{ij}(\thnn)+\phi_{1ij}(b^*_{i1})a_{n_i}^{-1}t,
$
where as before, $b^*_{i1}=\thnn+ca_{n_i}^{-1}t$ for some $0<c<1$. Accordingly we obtain, 
\beqn
\tau_{n_i}^{-1}\underset{|t|\le K\tau_{n_i}}{\sup}|S^*_{n_i}(t)|&=&\tau_{n_i}^{-1}\underset{|t|\le K\tau_{n_i}}{\sup} \left|a_{n_i}^{-1}\jsum w_{ij} \phi_{1ij}(b^*_{i1})a_{n_i}^{-1}t-\frac{t}{\sigma^2} \right|\\
&=&K\left| a_{n_i}^{-2}\underset{|t|\le K\tau_{n_i}}{\sup}\jsum w_{ij}\phi_{1ij}(b^*_{i1})-\frac{1}{\sigma^2}\right|\\
&\le&K\left| \tau_{n_i}a_{n_i}^{-2}\underset{|t|\le K\tau_{n_i}}{\sup}\jsum W_{ij}\phi_{1ij}(b^*_{i1})\right|+ K\left|a_{n_i}^{-2}\underset{|t|\le K\tau_{n_i}}{\sup}\jsum \phi_{1ij}(b^*_{i1})-\frac{1}{\sigma^2}\right|.
\eeqn
Further,
\beqn
E^*\left[\tau_{n_i}^{-1}\underset{|t|\le K\tau_{n_i}}{\sup}|S^*_{n_i}(t)|\right]^2&\le&K^2E^*\left| \tau_{n_i}a_{n_i}^{-2}\underset{|t|\le K\tau_{n_i}}{\sup}\jsum W_{ij}\phi_{1ij}(b^*_{i1})\right|^2\\
&+& K^2E^*\left|a_{n_i}^{-2}\underset{|t|\le K\tau_{n_i}}{\sup}\jsum \phi_{1ij}(b^*_{i1})-\frac{1}{\sigma^2}\right|^2\\
&+&K^2E^*\left| \tau_{n_i}a_{n_i}^{-2}\underset{|t|\le K\tau_{n_i}}{\sup}\jsum W_{ij}\phi_{1ij}(b^*_{i1})\right|\left|a_{n_i}^{-2}\underset{|t|\le K\tau_{n_i}}{\sup}\jsum \phi_{1ij}(b^*_{i1})-\frac{1}{\sigma^2}\right|.
\eeqn
By Lemma \ref{7.4.2} and Lemma \ref{7.3.1}, we have 
\[
E^*\left| \tau_{n_i}a_{n_i}^{-2}\underset{|t|\le K\tau_{n_i}}{\sup}\jsum W_{ij}\phi_{1ij}(b^*_{i1})\right|^2\to0,
\]
and
\[
E^*\left|a_{n_i}^{-2}\underset{|t|\le K\tau_{n_i}}{\sup}\jsum \phi_{1ij}(b^*_{i1})-\frac{1}{\sigma^2}\right|^2\to0.
\]
Thus, by an application of the Cauchy-Schwarz inequality we have proved (\ref{3.3}). 
Next, in light of (\ref{7.12}) we define
\beqn
A^*_{n_i}(t):=a_{n_i}^{-1}t\jsum w_{ij}\phi_{ij}(\thnn+a_{n_i}^{-1}t)= tS^*_{n_i}(t)+a_{n_i}^{-1}t\jsum w_{ij}\phi_{ij}(\thnn)+\frac{t^2}{\sigma^2}.  
\eeqn
Accordingly,
\beqn
\underset{|t|=K\tau_{n_i}}{\inf}A^*_{n_i}(t)\ge -K\tau_{n_i}\underset{|t|=K\tau_{n_i}}{\sup}|S^*_{n_i}(t)|-K\tau_{n_i}a_{n_i}^{-1}\left|\jsum w_{ij} \phi_{ij}(\thnn)\right|+\frac{K^2\tau_{n_i}^2}{\sigma^2}.
\eeqn
Recall that by Lemma \ref{7.4.3}, there exists a $K>0$ such that for any $\e>0$, 
\be\label{7.11}
P^*\l\left|a_{n_i}^{-1}\jsum W_{ij}\phi_{ij}(\thnn)\right|>K\r<\frac{\e}{2}.
\ee
Accordingly, by (\ref{7.11}) and (\ref{3.3}) we may choose  large enough $K$ such that for sufficiently large $n_i$,
\beqn
P^*\left(\underset{|t|=K\tau_{n_i}}{\inf}A_{n_i}(t)\ge0\right)&\ge& P^*\l\underset{|t|=K\tau_{n_i}}{\sup}|S^*_{n_i}(t)|+a_{n_i}^{-1}\left|\jsum w_{ij}\phi_{ij}(\thnn)\right|\le \frac{K\tau_{n_i}}{\sigma^2}\r\\
&=&P^*\l\underset{|t|=K\tau_{n_i}}{\sup}|S^*_{n_i}(t)|+a_{n_i}^{-1}\tau_{n_i}\left|\jsum W_{ij}\phi_{ij}(\thnn)\right|\le \frac{K\tau_{n_i}}{\sigma^2}\r\\
&=&1-P^*\l\underset{|t|=K\tau_{n_i}}{\sup}|S^*_{n_i}(t)|+a_{n_i}^{-1}\tau_{n_i}\left|\jsum W_{ij}\phi_{ij}(\thnn)\right|> \frac{K\tau_{n_i}}{\sigma^2}\r\\
&\ge&1-P^*\l\tau_{n_i}^{-1}\underset{|t|=K\tau_{n_i}}{\sup}|S^*_{n_i}(t)|>\frac{K}{4\sigma^2}\r-P^*\l a_{n_i}^{-1}\left|\jsum W_{ij}\phi_{ij}(\thnn)\right|>\frac{K}{4\sigma^2}\r\\
&\ge& 1-\e.
\eeqn
From  the continuity of $\jsum \phi_{ij}(\thtt)$ in $\thtt$, we have for sufficiently large $n_i$, that there exists a $K$ such that the equation 
$
\jsum w_{ij}\phi_{ij}(\thnn+a_{n_i}^{-1}t)=0, 
$
has a root, $t=T^*_{ni}$ in $|t|\le K\tau_{n_i}$, with a probability larger than $1-\e$. That is, we have
\[
\thss=\thnn+a_{ni}^{-1}T^*_{ni}, 
\]
where $|\tau_{n_i}^{-1}T^*_{ni}|<K$ in probability. Accordingly we may rewrite $\hat\theta^*_{_{STS}}$ as, 
\beqn
\hat\theta^*_{_{STS}}&=&\frac{1}{N}\isum u_i\thnn +\frac{1}{N}\isum u_ia_{ni}^{-1}T^*_{ni}\\
&=&\frac{1}{N}\isum u_i(\thzz+b_i+a_{ni}^{-1}T_{ni})+\frac{1}{N}\isum u_ia_{ni}^{-1}T^*_{ni}\\
&=&\frac{1}{N}\isum u_i\thzz+\frac{1}{N}\isum u_ib_i+\frac{1}{N}\isum u_ia_{ni}^{-1}T_{ni}+\frac{1}{N}\isum u_ia_{ni}^{-1}T^*_{ni}.
\eeqn
That is,
\[
\hat\theta^*_{_{STS}}-\thzz=\frac{1}{N}\isum (u_i-1)\thzz+\frac{1}{N}\isum u_ib_i+\frac{1}{N}\isum u_ia_{ni}^{-1}T_{ni}+\frac{1}{N}\isum u_ia_{ni}^{-1}T^*_{ni}.
\]
Additionally, by Lemma \ref{2.4.6},  we have $\frac{1}{N}\isum (u_i-1)\overset{p^*}{\to}0, $ as well as, $\frac{1}{N}\isum u_ib_i\overset{p^*}{\to}0$. Further, we also have that 
\beqn
\frac{1}{N}\isum u_ia_{ni}^{-1}T_{ni}= \frac{1}{N}\isum (u_i-1)a_{ni}^{-1}T_{ni}+\frac{1}{N}\isum a_{ni}^{-1}T_{ni}.
\eeqn
Now by Lemma \ref{7.3.2}  and the fact $T_{ni}=O_p(1)$, we obtain, with $U_i:=(u_i-1)/\tau_N$,  that 
\beqn
E^*(\frac{1}{N}\isum (u_i-1)a_{ni}^{-1}T_{ni})^2&=&E^*(\frac{\tau_N}{N}\isum U_ia_{ni}^{-1}T_{ni})^2\\
&\le&\frac{\tau_N^2}{N^2}\isum a_{ni}^{-2}T^2_{ni}+(N-1)O(\frac{1}{N})\frac{\tau_N^2}{N^2}\isum a_{ni}^{-2}T^2_{ni} \overset{p}{\to}0,
\eeqn
as well as, $\frac{1}{N}\isum a_{ni}^{-1}T_{ni}\overset{p}{\to} 0$. That is, we have established that, $E^*(\frac{1}{N}\isum u_ia_{ni}^{-1}T_{ni})^2\overset{p}{\to} 0$. Accordingly we conclude, $P^*(|\frac{1}{N}\isum u_ia_{ni}^{-1}T_{ni}|>\e)=o_p(1)$. 
Similarly,
\beqn
\frac{1}{N}\isum u_ia_{ni}^{-1}T^*_{ni}= \frac{1}{N}\isum (u_i-1)a_{ni}^{-1}T^*_{ni}+\frac{1}{N}\isum a_{ni}^{-1}T^*_{ni},
\eeqn
where by Lemma \ref{7.3.2}, {\it Assumption B (3)}  and the fact $\tau_{n_i}^{-1}T^{*}_{ni}=O_{p^*}(1)$, we obtain, 
\beqn
E^*(\frac{1}{N}\isum (u_i-1)a_{ni}^{-1}T^*_{ni})^2&=&E^*(\frac{\tau_N}{N}\isum U_ia_{ni}^{-1}T^*_{ni})^2\\
&\le&\frac{\tau_N^2}{N^2}\isum a_{ni}^{-2}T^{*2}_{ni}+(N-1)O(\frac{1}{N})\frac{\tau_N^2}{N^2}\isum a_{ni}^{-2}T^{*2}_{ni}\\
&=&(1+(N-1)O(\frac{1}{N}))\frac{\tau_N^2}{N^2}\isum \frac{\tau^2_{n_i}}{a_{ni}^{2}}\tau_{n_i}^{-2}T^{*2}_{ni}\overset{p}{\to}0.
\eeqn
Finally, by Lemma \ref{7.3.2},
\beqn
\frac{1}{N}\isum a_{ni}^{-1}T^*_{ni}&=& \frac{1}{N}\isum \frac{\tau_{n_i}}{a_{ni}}\tau_{n_i}^{-1}T^*_{ni}\to 0.
\eeqn
Accordingly we also conclude that, $P^*(|\frac{1}{N}\isum u_ia_{ni}^{-1}T^*_{ni}|>\e)=o_p(1)$. Hence, we have proved that  $P^*(|\hat\theta^*_{_{STS}}-\thzz|>\e)=o_p(1)$. 
\qed
\bigskip

\noindent For the related asymptotic normality results as stated in Theorem \ref{5.1.2}, we need the following two Lemmas.
\begin{lemma}\label{7.4.4}
	Suppose that the conditions of {\it Assumptions A } and {\it B} hold. If $\frac{\tau_{n_i}}{\tau_N}=o(\sqrt{n_i})$ then as $n_i\to\infty$ and $N\to\infty$,
	\[
	\frac{\tau_N^{-1}}{\sqrt{N}}\isum u_ia_{n_i}^{-2}\jsum w_{ij}\phi_{ij}(\thnn)\overset{p^*}{\to} 0.
	\]
\end{lemma}
{\it \underbar{Proof of Lemma \ref{7.4.4}}:} Let
\[
X^*_{ni}:= \tau_{n_i}^{-1}a_{n_i}^{-1}\jsum w_{ij} \phi_{ij}(\thnn)=a_{n_i}^{-1}\jsum W_{ij} \phi_{ij}(\thnn).
\]
Clearly $E^*(X^*_{ni})=0$, and $X^*_{n_i}$ are independent for $i$ in $1,2,\dots, N$. Further, by Lemma \ref{7.4.1} we have, as $n_i\to\infty$, that 
\beqn
E^*(X^{*2}_{ni})&=& E^*(a_{n_i}^{-1}\jsum W_{ij} \phi_{ij}(\thnn))^2\\
&=&a_{n_i}^{-2}\l\jsum \phi^2_{ij}(\thnn)+O(\frac{1}{n_i})\underset{j_1\ne j_2}{\sum}\phi_{ij_1}(\thnn)\phi_{ij_2}(\thnn)\r\\
&=&a_{n_i}^{-2}\l\jsum \phi^2_{ij}(\thnn)+O(\frac{1}{n_i})\ll\jsum \phi_{ij}(\thnn)\rr^2-O(\frac{1}{n_i})\jsum \phi^2_{ij}(\thnn)\r\\
&=& (1-O(\frac{1}{n_i}))a_{n_i}^{-2}\jsum \phi^2_{ij}(\thnn)\to1.
\eeqn
Thus, with $U_i=(u_i-1)/\sqrt{\tau_N}$, 
\beqn
\frac{\tau_N^{-1}}{\sqrt{N}}\isum u_ia_{n_i}^{-2}\jsum w_{ij}\phi_{ij}(\thnn)&=& \frac{\tau_N^{-1}}{\sqrt{N}}\isum u_ia_{n_i}^{-1}\tau_{n_i}X^*_{ni}\\
&=&\frac{1}{\sqrt{N}}\isum U_ia_{n_i}^{-1}\tau_{n_i}X^*_{ni}+\frac{\tau_N^{-1}}{\sqrt{N}}\isum a_{n_i}^{-1}\tau_{n_i}X^*_{ni}.
\eeqn
Since $U_i$ and $X^*_{ni}$ are independent, we obtain, 
\beqn
E^*(\frac{1}{\sqrt{N}}\isum U_ia_{n_i}^{-1}\tau_{n_i}X^*_{ni})^2&=&\frac{1}{N}\isum E^* (U_i^2a_{n_i}^{-2}\tau_{n_i}^2X^{*2}_{ni})\\
&+&\underset{i_1\ne i_2}{\sum}E^*(U_{i_1}U_{i_2}a_{n_{i_1}}^{-1}a_{n_{i_2}}^{-1}\tau_{n_{i_1}}\tau_{n_{i_2}}X^{*}_{ni_1}X^{*}_{ni_2})\\
&=&\frac{1}{N}\isum a_{n_i}^{-2}\tau_{n_i}^2E^* (X^{*2}_{ni})\to 0.
\eeqn
Finally, since  $\frac{\tau_{n_i}}{\tau_N}=o(\sqrt{n_i})$, we also have, 
\beqn
E^*(\frac{\tau_N^{-1}}{\sqrt{N}}\isum a_{n_i}^{-1}\tau_{n_i}X^*_{ni})^2&=&\frac{\tau_N^{-2}}{N}\isum a_{n_i}^{-2}\tau_{n_i}^2E^*(X^{*2}_{ni})\\
&=&\frac{1}{N}\isum \frac{\tau_{n_i}^2}{\tau_N^{2}}a_{n_i}^{-2}E^*(X^{*2}_{ni})\to 0.
\eeqn
Accordingly we obtain that, 
\[
\frac{\tau_N^{-1}}{\sqrt{N}}\isum u_ia_{n_i}^{-2}\jsum w_{ij}\phi_{ij}(\thnn)\overset{p^*}{\to} 0.
\]
\qed
\begin{lemma}\label{7.4.5}
	Suppose that the conditions of {\it Assumptions A } and {\it B} hold.  If $\frac{\tau_{n_i}}{\tau_N}=o(\sqrt{n_i})$ then as $n_i\to\infty$ and $N\to\infty$,
	\[
	\frac{\lambda^{-1}\tau_N^{-1}\sigma^2}{\sqrt{N}}\isum u_ia_{n_i}^{-1}S_{n_i}(T^*_{ni}) \overset{p^*}{\to} 0.
	\]
\end{lemma}
{\it \underbar{Proof of Lemma \ref{7.4.5}}:} We first write
\beqn
\frac{\lambda^{-1}\tau_N^{-1}\sigma^2}{\sqrt{N}}\isum u_ia_{n_i}^{-1}S_{n_i}(T^*_{ni})=\frac{\lambda^{-1}\sigma^2}{\sqrt{N}}\isum U_ia_{n_i}^{-1}S_{n_i}(T^*_{ni})+\frac{\lambda^{-1}\tau_N^{-1}\sigma^2}{\sqrt{N}}\isum a_{n_i}^{-1}S_{n_i}(T^*_{ni}).
\eeqn
By Lemma \ref{7.3.2}, {\it Assumption B (3)} and the fact $\tau_N^{-1}S_{n_i}(T^*_{ni})\overset{p^*}{\to} 0$,
\[
\frac{\lambda^{-1}\tau_N^{-1}\sigma^2}{\sqrt{N}}\isum a_{n_i}^{-1}S_{n_i}(T^*_{ni})\overset{p^*}{\to} 0.
\]
Further, it can be seen that, 
\beqn
E^*(\frac{1}{\sqrt{N}}\isum U_ia_{ni}^{-1}S_{n_i}(T^*_{ni}))^2&\le&\frac{1}{N}\l1+(N-1)O(\frac{1}{N})\r\isum a_{ni}^{-2}E^*(S^2_{n_i}(T^*_{ni}))\to 0. 
\eeqn
Thus we have,
\[
\frac{\lambda^{-1}\tau_N^{-1}\sigma^2}{\sqrt{N}}\isum u_ia_{n_i}^{-1}S_{n_i}(T^*_{ni}) \overset{p^*}{\to} 0.
\]
\qed

\noindent{\underline{\it Proof of Theorem \ref{5.1.2}:} } By Theorem \ref{5.1.1} and (\ref{7.12}) we express,
\[
\thss-\thnn=-a_{ni}^{-2}\sigma^2\jsum w_{ij}\phi_{ij}(\thnn)- a_{n_i}^{-1}\sigma^2S_{n_i}(T^*_{ni}). 
\]
Accordingly we have,
\[
\hat\theta^*_{_{STS}}-\hat\theta_{_{STS}}= \frac{1}{N}\isum (u_i-1)\thnn -\frac{\sigma^2}{N}\isum u_ia_{n_i}^{-2}\jsum w_{ij}\phi_{ij}(\thnn)-\frac{\sigma^2}{N}\isum u_ia_{n_i}^{-1}S_{n_i}(T^*_{ni}), 
\]
where $|T^*_{ni}|<K\tau_{n_i}$ in probability. Further,
\beqn
\lambda^{-1}\tau_N^{-1}\sqrt{N}(\hat\theta^*_{_{STS}}-\hat\theta_{_{STS}})&=& \frac{\lambda^{-1}\tau_N^{-1}}{\sqrt{N}}\isum (u_i-1)\thnn\\
&-&\frac{\lambda^{-1}\tau_N^{-1}\sigma^2}{\sqrt{N}}\isum u_ia_{n_i}^{-2}\jsum w_{ij}\phi_{ij}(\thnn)\\
&-&\frac{\lambda^{-1}\tau_N^{-1}\sigma^2}{\sqrt{N}}\isum u_ia_{n_i}^{-1}S_{n_i}(T^*_{ni})\\
&\equiv& I_1+I_2+I_3 .
\eeqn
By Lemma \ref{7.4.4},  $I_2\overset{p^*}{\to} 0$, and by Lemma \ref{7.4.5}, $I_3\overset{p^*}{\to} 0$,  and  therefore  it remains only to consider $I_1$. Now, observe that, 
\beqn
I_1:= \frac{\lambda^{-1}\tau_N^{-1}}{\sqrt{N}}\isum (u_i-1)\thnn=\frac{\lambda^{-1}}{\sqrt{N}}\isum U_i(b_i+\thzz)+\frac{\lambda^{-1}}{\sqrt{N}}\isum U_ia_{n_i}^{-1}T_{ni}.
\eeqn
By Lemma \ref{7.3.2},
\beqn
E^*(\frac{1}{\sqrt{N}}\isum U_ia_{ni}^{-1}T_{ni})^2&\le&\frac{1}{N}\isum a_{ni}^{-2}T^2_{ni}+(N-1)O(\frac{1}{N})\frac{1}{N}\isum a_{ni}^{-2}T^2_{ni} \overset{p}{\to}0.
\eeqn
Further by Lemma \ref{2.4.3}, 
\[
\bar{U}_N:=\frac{1}{N}\isum U_i\equiv \frac{1}{N} \isum \frac{u_i-1}{\tau_N}\overset{p^*}{\to} 0,
\]
and clearly, $\sqrt{N}(\bar{b}+\thzz)\Rightarrow {\Cal{N}}(\thzz,\lambda^2)$.  Accordingly we have, $\frac{\lambda^{-1}}{N}\isum(b_i-\bar{b})^2\to1 \ \ a.s.$ as well as $\sqrt{N} \bar{U}(\bar{b}+\thzz)\overset{p^*}{\to} 0$. Further, by Lemma 4.6 of Praestgaard and Wellner (1993), we have that 
\[
\frac{\lambda^{-1}}{\sqrt{N}}\isum U_i(b_i+\thzz)\Rightarrow {\cal N}(0,1).
\]
Thus we have 
\[
\frac{\lambda^{-1}\tau_N^{-1}}{\sqrt{N}}\isum (u_i-1)\thnn\Rightarrow {\Cal{N}}(0,1).
\]
Finally we conclude that as $n_i\to\infty$ and $N\to \infty$, 
\[
\lambda^{-1}\tau_N^{-1}\sqrt{N}(\hat\theta^*_{_{STS}}-\hat\theta_{_{STS}})\Rightarrow {\Cal{N}}(0,1).
\]
\qed

\vfill\eject

\end{document}